\documentclass[12pt,reqno]{amsart}

\usepackage{amssymb}

\catcode`\@=11

\thinlines
\newskip\Einheit \Einheit=0.6cm
\newcount\xcoord \newcount\ycoord
\newdimen\xdim \newdimen\ydim \newdimen\PfadD@cke \newdimen\Pfadd@cke
\def\PfadDicke#1{\PfadD@cke#1 \divide\PfadD@cke by2 \Pfadd@cke\PfadD@cke \multiply\PfadD@cke by2}
\long\def\LOOP#1\REPEAT{\def\BODY{#1}\ITERATE}
\def\ITERATE{\BODY \let\next\ITERATE \else\let\next\relax\fi \next}
\let\REPEAT=\fi
\def\Punkt{\hbox{\raise-2pt\hbox to0pt{\hss\scriptsize$\bullet$\hss}}}
\def\DuennPunkt(#1,#2){\unskip
  \raise#2 \Einheit\hbox to0pt{\hskip#1 \Einheit
          \raise-2.5pt\hbox to0pt{\hss\normalsize$\bullet$\hss}\hss}}
\def\NormalPunkt(#1,#2){\unskip
  \raise#2 \Einheit\hbox to0pt{\hskip#1 \Einheit
          \raise-3pt\hbox to0pt{\hss\large$\bullet$\hss}\hss}}
\def\DickPunkt(#1,#2){\unskip
  \raise#2 \Einheit\hbox to0pt{\hskip#1 \Einheit
          \raise-4pt\hbox to0pt{\hss\Large$\bullet$\hss}\hss}}
\def\Kreis(#1,#2){\unskip
  \raise#2 \Einheit\hbox to0pt{\hskip#1 \Einheit
          \raise-4pt\hbox to0pt{\hss\Large$\circ$\hss}\hss}}
\def\Diagonale(#1,#2)#3{\unskip\leavevmode
  \xcoord#1\relax \ycoord#2\relax
      \raise\ycoord \Einheit\hbox to0pt{\hskip\xcoord \Einheit
         \unitlength\Einheit
         \line(1,1){#3}\hss}}
\def\AntiDiagonale(#1,#2)#3{\unskip\leavevmode
  \xcoord#1\relax \ycoord#2\relax 
      \raise\ycoord \Einheit\hbox to0pt{\hskip\xcoord \Einheit
         \unitlength\Einheit
         \line(1,-1){#3}\hss}}
\def\Pfad(#1,#2),#3\endPfad{\unskip\leavevmode
  \xcoord#1 \ycoord#2 \thicklines\ZeichnePfad#3\endPfad\thinlines}
\def\ZeichnePfad#1{\ifx#1\endPfad\let\next\relax
  \else\let\next\ZeichnePfad
    \ifnum#1=1
      \raise\ycoord \Einheit\hbox to0pt{\hskip\xcoord \Einheit
         \vrule height\Pfadd@cke width1 \Einheit depth\Pfadd@cke\hss}%
      \advance\xcoord by 1
    \else\ifnum#1=2
      \raise\ycoord \Einheit\hbox to0pt{\hskip\xcoord \Einheit
        \hbox{\hskip-\PfadD@cke\vrule height1 \Einheit width\PfadD@cke depth0pt}\hss}%
      \advance\ycoord by 1
    \else\ifnum#1=3
      \raise\ycoord \Einheit\hbox to0pt{\hskip\xcoord \Einheit
         \unitlength\Einheit
         \line(1,1){1}\hss}
      \advance\xcoord by 1
      \advance\ycoord by 1
    \else\ifnum#1=4
      \raise\ycoord \Einheit\hbox to0pt{\hskip\xcoord \Einheit
         \unitlength\Einheit
         \line(1,-1){1}\hss}
      \advance\xcoord by 1
      \advance\ycoord by -1
    \else\ifnum#1=5
      \advance\xcoord by -1
      \raise\ycoord \Einheit\hbox to0pt{\hskip\xcoord \Einheit
         \vrule height\Pfadd@cke width1 \Einheit depth\Pfadd@cke\hss}%
    \else\ifnum#1=6
      \advance\ycoord by -1
      \raise\ycoord \Einheit\hbox to0pt{\hskip\xcoord \Einheit
        \hbox{\hskip-\PfadD@cke\vrule height1 \Einheit width\PfadD@cke depth0pt}\hss}%
    \else\ifnum#1=7
      \advance\xcoord by -1
      \advance\ycoord by -1
      \raise\ycoord \Einheit\hbox to0pt{\hskip\xcoord \Einheit
         \unitlength\Einheit
         \line(1,1){1}\hss}
    \else\ifnum#1=8
      \advance\xcoord by -1
      \advance\ycoord by +1
      \raise\ycoord \Einheit\hbox to0pt{\hskip\xcoord \Einheit
         \unitlength\Einheit
         \line(1,-1){1}\hss}
    \fi\fi\fi\fi
    \fi\fi\fi\fi
  \fi\next}
\def\hSSchritt{\leavevmode\raise-.4pt\hbox to0pt{\hss.\hss}\hskip.2\Einheit
  \raise-.4pt\hbox to0pt{\hss.\hss}\hskip.2\Einheit
  \raise-.4pt\hbox to0pt{\hss.\hss}\hskip.2\Einheit
  \raise-.4pt\hbox to0pt{\hss.\hss}\hskip.2\Einheit
  \raise-.4pt\hbox to0pt{\hss.\hss}\hskip.2\Einheit}
\def\vSSchritt{\vbox{\baselineskip.2\Einheit\lineskiplimit0pt
\hbox{.}\hbox{.}\hbox{.}\hbox{.}\hbox{.}}}
\def\DSSchritt{\leavevmode\raise-.4pt\hbox to0pt{%
  \hbox to0pt{\hss.\hss}\hskip.2\Einheit
  \raise.2\Einheit\hbox to0pt{\hss.\hss}\hskip.2\Einheit
  \raise.4\Einheit\hbox to0pt{\hss.\hss}\hskip.2\Einheit
  \raise.6\Einheit\hbox to0pt{\hss.\hss}\hskip.2\Einheit
  \raise.8\Einheit\hbox to0pt{\hss.\hss}\hss}}
\def\dSSchritt{\leavevmode\raise-.4pt\hbox to0pt{%
  \hbox to0pt{\hss.\hss}\hskip.2\Einheit
  \raise-.2\Einheit\hbox to0pt{\hss.\hss}\hskip.2\Einheit
  \raise-.4\Einheit\hbox to0pt{\hss.\hss}\hskip.2\Einheit
  \raise-.6\Einheit\hbox to0pt{\hss.\hss}\hskip.2\Einheit
  \raise-.8\Einheit\hbox to0pt{\hss.\hss}\hss}}
\def\SPfad(#1,#2),#3\endSPfad{\unskip\leavevmode
  \xcoord#1 \ycoord#2 \ZeichneSPfad#3\endSPfad}
\def\ZeichneSPfad#1{\ifx#1\endSPfad\let\next\relax
  \else\let\next\ZeichneSPfad
    \ifnum#1=1
      \raise\ycoord \Einheit\hbox to0pt{\hskip\xcoord \Einheit
         \hSSchritt\hss}%
      \advance\xcoord by 1
    \else\ifnum#1=2
      \raise\ycoord \Einheit\hbox to0pt{\hskip\xcoord \Einheit
        \hbox{\hskip-2pt \vSSchritt}\hss}%
      \advance\ycoord by 1
    \else\ifnum#1=3
      \raise\ycoord \Einheit\hbox to0pt{\hskip\xcoord \Einheit
         \DSSchritt\hss}
      \advance\xcoord by 1
      \advance\ycoord by 1
    \else\ifnum#1=4
      \raise\ycoord \Einheit\hbox to0pt{\hskip\xcoord \Einheit
         \dSSchritt\hss}
      \advance\xcoord by 1
      \advance\ycoord by -1
    \else\ifnum#1=5
      \advance\xcoord by -1
      \raise\ycoord \Einheit\hbox to0pt{\hskip\xcoord \Einheit
         \hSSchritt\hss}%
    \else\ifnum#1=6
      \advance\ycoord by -1
      \raise\ycoord \Einheit\hbox to0pt{\hskip\xcoord \Einheit
        \hbox{\hskip-2pt \vSSchritt}\hss}%
    \else\ifnum#1=7
      \advance\xcoord by -1
      \advance\ycoord by -1
      \raise\ycoord \Einheit\hbox to0pt{\hskip\xcoord \Einheit
         \DSSchritt\hss}
    \else\ifnum#1=8
      \advance\xcoord by -1
      \advance\ycoord by 1
      \raise\ycoord \Einheit\hbox to0pt{\hskip\xcoord \Einheit
         \dSSchritt\hss}
    \fi\fi\fi\fi
    \fi\fi\fi\fi
  \fi\next}
\def\Koordinatenachsen(#1,#2){\unskip
 \hbox to0pt{\hskip-.5pt\vrule height#2 \Einheit width.5pt depth1 \Einheit}%
 \hbox to0pt{\hskip-1 \Einheit \xcoord#1 \advance\xcoord by1
    \vrule height0.25pt width\xcoord \Einheit depth0.25pt\hss}}
\def\Koordinatenachsen(#1,#2)(#3,#4){\unskip
 \hbox to0pt{\hskip-.5pt \ycoord-#4 \advance\ycoord by1
    \vrule height#2 \Einheit width.5pt depth\ycoord \Einheit}%
 \hbox to0pt{\hskip-1 \Einheit \hskip#3\Einheit 
    \xcoord#1 \advance\xcoord by1 \advance\xcoord by-#3 
    \vrule height0.25pt width\xcoord \Einheit depth0.25pt\hss}}
\def\Gitter(#1,#2){\unskip \xcoord0 \ycoord0 \leavevmode
  \LOOP\ifnum\ycoord<#2
    \loop\ifnum\xcoord<#1
      \raise\ycoord \Einheit\hbox to0pt{\hskip\xcoord \Einheit\Punkt\hss}%
      \advance\xcoord by1
    \repeat
    \xcoord0
    \advance\ycoord by1
  \REPEAT}
\def\Gitter(#1,#2)(#3,#4){\unskip \xcoord#3 \ycoord#4 \leavevmode
  \LOOP\ifnum\ycoord<#2
    \loop\ifnum\xcoord<#1
      \raise\ycoord \Einheit\hbox to0pt{\hskip\xcoord \Einheit\Punkt\hss}%
      \advance\xcoord by1
    \repeat
    \xcoord#3
    \advance\ycoord by1
  \REPEAT}
\def\Label#1#2(#3,#4){\unskip \xdim#3 \Einheit \ydim#4 \Einheit
  \def\lo{\advance\xdim by-.5 \Einheit \advance\ydim by.5 \Einheit}%
  \def\llo{\advance\xdim by-.25cm \advance\ydim by.5 \Einheit}%
  \def\loo{\advance\xdim by-.5 \Einheit \advance\ydim by.25cm}%
  \def\o{\advance\ydim by.25cm}%
  \def\ro{\advance\xdim by.5 \Einheit \advance\ydim by.5 \Einheit}%
  \def\rro{\advance\xdim by.25cm \advance\ydim by.5 \Einheit}%
  \def\roo{\advance\xdim by.5 \Einheit \advance\ydim by.25cm}%
  \def\l{\advance\xdim by-.30cm}%
  \def\r{\advance\xdim by.30cm}%
  \def\lu{\advance\xdim by-.5 \Einheit \advance\ydim by-.6 \Einheit}%
  \def\llu{\advance\xdim by-.25cm \advance\ydim by-.6 \Einheit}%
  \def\luu{\advance\xdim by-.5 \Einheit \advance\ydim by-.30cm}%
  \def\u{\advance\ydim by-.30cm}%
  \def\ru{\advance\xdim by.5 \Einheit \advance\ydim by-.6 \Einheit}%
  \def\rru{\advance\xdim by.25cm \advance\ydim by-.6 \Einheit}%
  \def\ruu{\advance\xdim by.5 \Einheit \advance\ydim by-.30cm}%
  #1\raise\ydim\hbox to0pt{\hskip\xdim
     \vbox to0pt{\vss\hbox to0pt{\hss$#2$\hss}\vss}\hss}%
}
\catcode`\@=12

\usepackage{color}

\input texdraw

\def\Ringerl(#1 #2){\move(#1 #2)\fcir f:0 r:.15}
\def\ringerl(#1 #2){\move(#1 #2)\fcir f:0 r:.1}
\def\Mark(#1 #2){\move(#1 #2)\fcir f:0 r:.2}

\newtheorem{theorem}{Theorem}
\newtheorem{proposition}[theorem]{Proposition}
\newtheorem{lemma}[theorem]{Lemma}
\newtheorem{corollary}[theorem]{Corollary}

\theoremstyle{remark}
\newtheorem{remark}{Remark}

\numberwithin{equation}{section}

\newcounter{saveeqn}
\newcommand{\alphaeqn}{\setcounter{saveeqn}{\value{equation}}%
\setcounter{equation}{0}%
\global\def\theequation{\mbox{\thesection.\arabic{saveeqn}\alph{equation}}}}
\newcommand{\reseteqn}{\setcounter{equation}{\value{saveeqn}}%
\global\def\theequation{\thesection.\arabic{equation}}}

\def\({\left(}
\def\){\right)}

\def\Fix{\operatorname{Fix}}
\def\Cent{\operatorname{Cent}}
\def\Cat{\operatorname{Cat}}

\def\rk{\operatorname{rk}}

\def\ep{\varepsilon}
\def\al{\alpha}
\def\be{\beta}
\def\ga{\gamma}

\begin{document}

\title[Cyclic sieving for generalised non-crossing partitions]
{Cyclic sieving for generalised non-crossing partitions associated
with complex reflection groups of exceptional type}
\newbox\Aut
\setbox\Aut\vbox{
\centerline{\sc 
Christian Krattenthaler$^{\dagger}$ \rm and \sc Thomas W.
M\"uller$^\ddagger$}
\vskip18pt
\centerline{$^\dagger$ Fakult\"at f\"ur Mathematik, Universit\"at Wien,}
\centerline{Nordbergstra\ss e 15, A-1090 Vienna, Austria.}
\centerline{\footnotesize WWW: \footnotesize\tt
http://www.mat.univie.ac.at/\lower0.5ex\hbox{\~{}}kratt} 
\vskip18pt
\centerline{$^\ddagger$ School of Mathematical Sciences,}
\centerline{Queen Mary \& Westfield College, University of London,}
\centerline{Mile End Road, London E1 4NS, United Kingdom.}
\centerline{\footnotesize WWW: \tt http://www.maths.qmw.ac.uk/\lower0.5ex\hbox{\~{}}twm/}
}
\author[C. Krattenthaler and T. W. M\"uller]{\box\Aut}

\address{Fakult\"at f\"ur Mathematik, Universit\"at Wien,
Nordbergstra{\ss}e~15, A-1090 Vienna, Austria.
WWW: \tt http://www.mat.univie.ac.at/\lower0.5ex\hbox{\~{}}kratt.}

\address{School of Mathematical Sciences, Queen Mary
\& Westfield College, University of London,
Mile End Road, London E1 4NS, United Kingdom.\newline
\tt http://www.maths.qmw.ac.uk/\lower0.5ex\hbox{\~{}}twm/.}


\dedicatory{Dedicated to the memory of Herb Wilf}


\thanks{$^\dagger$Research partially supported by the Austrian
Science Foundation FWF, grants Z130-N13 and S9607-N13,
the latter in the framework of the National Research Network
``Analytic Combinatorics and Probabilistic Number Theory."\newline
\indent
$^\ddagger$Research supported by the Austrian
Science Foundation FWF, Lise Meitner grant M1201-N13}

\begin{abstract}
We prove that the generalised non-crossing partitions associated with
well-generated complex reflection groups of exceptional type
obey two different cyclic sieving phenomena, as conjectured by
Armstrong, and by Bessis and Reiner.
The computational details are provided in the manuscript
{\it``Cyclic sieving for generalised non-crossing partitions associated
with complex reflection groups of exceptional type --- the details"}
[{\tt ar$\chi$iv:1001.0030}].
\end{abstract}

\maketitle

\section{Introduction}

In his memoir \cite{ArmDAA}, Armstrong introduced {\em generalised
non-crossing partitions} associated with finite (real) reflection groups,
thereby embedding Kreweras' non-crossing partitions \cite{KrewAC},
Edelman's $m$-divisible non-crossing partitions \cite{EdelAA}, 
the non-crossing partitions associated with reflection groups due to
Bessis \cite{BesDAA} and Brady and Watt \cite{BRWaAA} into one
uniform framework. Bessis and Reiner \cite{BeReAA} observed that
Armstrong's definition can be straightforwardly extended to {\em
well-generated complex reflection groups} (see Section~\ref{sec:prel} 
for the precise definition). 
These generalised non-crossing partitions possess a wealth of
beautiful properties, and they display deep and surprising 
relations to other combinatorial objects defined for reflection
groups (such as the generalised cluster complex of Fomin and
Reading \cite{FoReAA}, or the extended Shi arrangement and
the geometric multichains of filters
of Athanasiadis \cite{AthaAG,AthaAH}); 
see Armstrong's memoir \cite{ArmDAA} and the references given therein.

On the other hand,
{\it cyclic sieving} is a phenomenon brought to light by Reiner,
Stanton and White \cite{ReSWAA}.
It extends the so-called
``$(-1)$-phenomenon" of Stembridge \cite{StemAL,StemAP}. 
Cyclic sieving can be
defined in three equivalent ways (cf.\ \cite[Prop.~2.1]{ReSWAA}). 
The one which gives the name 
can be described as follows: given a set $S$ of combinatorial
objects, an action on $S$ of a cyclic group $G=\langle g\rangle$ with
generator $g$ of order $n$, and a polynomial $P(q)$ in $q$
with non-negative integer coefficients, we say
that the triple $(S,P,G)$ {\it exhibits the cyclic sieving
phenomenon}, if the number of elements of $S$ fixed by $g^k$ equals
$P(e^{2\pi ik/n})$. In \cite{ReSWAA} it is shown that this phenomenon
occurs in surprisingly many contexts, and several further instances
have been discovered since then.

In \cite[Conj.~5.4.7]{ArmDAA} (also appearing in \cite[Conj.~6.4]{BeReAA}) and 
\cite[Conj.~6.5]{BeReAA},
Armstrong, respectively Bessis and Reiner, conjecture that generalised
non-crossing partitions for irreducible 
well-generated complex reflection groups
exhibit two different cyclic sieving phenomena (see
Sections~\ref{sec:siev1} and \ref{sec:siev2} for the precise
statements). 

According to the classification of these groups due to
Shephard and Todd \cite{ShToAA}, there are two infinite families of
irreducible well-generated complex reflection groups, namely the
groups $G(d,1,n)$ and $G(e,e,n)$, where $n,d,e$ are positive integers, 
and there are 26 exceptional groups.
For the infinite families of types $G(d,1,n)$ and $G(e,e,n)$, the two
cyclic sieving conjectures follow from the results
in \cite{KratCG}. 

The purpose of the present article is to present a proof of the cyclic sieving
conjectures of Armstrong, and of Bessis and Reiner, for the
26 exceptional types, thus completing the proof of these conjectures. 
Since the generalised non-crossing partitions
feature a parameter $m$, from the outset 
this is {\it not\/} a finite problem. Consequently, we first need
several auxiliary results to reduce the conjectures for each of the
26 exceptional types to a {\it finite} problem. Subsequently, we use
Stembridge's {\sl Maple} package {\tt coxeter} \cite{StemAZ} and the
{\sl GAP} package {\tt CHEVIE} \cite{chevAA,MichAA} to carry out the remaining
{\it finite} computations. The details of these computations are provided
in \cite{KrMuAE}. In the present paper, we content ourselves with
exemplifying the necessary computations by going through some
representative cases.
It is interesting to observe that, for the
verification of the type $E_8$ case, it is essential to use the
decomposition numbers in the sense of \cite{KratCB,KratCF,KrMuAB} because,
otherwise, the necessary computations would not be feasible 
in reasonable time with the currently available computer facilities.
We point out that, for the special case where the aforementioned
parameter $m$ is equal to $1$, the first cyclic sieving conjecture
has been proved in a uniform fashion by Bessis and Reiner in \cite{BeReAA}. 
(See \cite{ArSTAA} for a --- non-uniform --- proof of cyclic sieving
for non-crossing partitions associated with {\it real\/} reflection groups
under the action of the so-called Kreweras map, a special case of
the second cyclic sieving phenomenon discussed in the present paper.)
The crucial result on which the proof of Bessis and Reiner 
is based is \eqref{eq:7}
below, and it plays an important role in our reduction of the
conjectures for the 26 exceptional groups to a finite problem.

Our paper is organised as follows.
In the next section, we recall the definition of generalised
non-crossing partitions for well-generated complex reflection
groups and of decomposition numbers in the sense of
\cite{KratCB,KratCF,KrMuAB}, and we review some basic facts. 
The first cyclic sieving conjecture is subsequently stated in
Section~\ref{sec:siev1}. In Section~\ref{sec:pol}, we outline an
elementary proof that the $q$-Fu\ss--Catalan number, which is
the polynomial $P$ in the cyclic sieving phenomena concerning
the generalised non-crossing partitions for well-generated 
complex reflection groups, is always a polynomial with non-negative
integer coefficients, as required by the definition of cyclic
sieving. (Full details can be found in \cite[Sec.~4]{KrMuAE}.
The reader is referred to the first paragraph of 
Section~\ref{sec:pol} for comments on other
approaches for establishing polynomiality with non-negative
coefficients.)
Section~\ref{sec:aux1} contains the
announced auxiliary results which, for the 26 exceptional types, 
allow a reduction of the conjecture to a finite problem. 
In Section~\ref{sec:Beweis1}, we discuss a few cases which,
in a representative manner, demonstrate
how to perform the
remaining case-by-case verification of the conjecture.
For full details, we refer the reader to \cite[Sec.~6]{KrMuAE}.
The second cyclic sieving conjecture is stated in
Section~\ref{sec:siev2}. Section~\ref{sec:aux2} contains the
auxiliary results which, for the 26 exceptional types, 
allow a reduction of the conjecture to a finite problem, while
in Section~\ref{sec:Beweis2} we discuss some representative cases of
the remaining case-by-case verification of the conjecture.
Again, for full details we refer the reader to \cite[Sec.~9]{KrMuAE}.

\section{Preliminaries}
\label{sec:prel}

A {\it complex reflection group} is a 
group generated by (complex) reflections in
$\mathbb C^n$. (Here, a reflection is a non-trivial element of
$GL_n(\mathbb C)$ which fixes a hyperplane pointwise and which 
has finite order.) We refer to \cite{LeTaAA} for an in-depth 
exposition of the theory complex reflection groups.

Shephard and Todd provided a complete classification of all {\it
finite} complex reflection groups in \cite{ShToAA} (see also
\cite[Ch.~8]{LeTaAA}). According to this
classification, an arbitrary complex reflection group $W$ decomposes into
a direct product of {\it irreducible} complex reflection groups, 
acting on mutually orthogonal subspaces of the complex vector space
on which $W$ is acting. Moreover, the list of irreducible complex
reflection groups consists of the infinite family of groups
$G(m,p,n)$, where $m,p,n$ are positive integers, and $34$ exceptional
groups, denoted $G_4,G_5,\dots,G_{37}$ by Shephard and Todd.

In this paper, we are only interested in finite complex
reflection groups which are {\it well-generated}. 
A complex reflection group of rank $n$ is called {\it well-generated\/} if
it is generated by $n$ reflections.\footnote{We refer to
\cite[Def.~1.29]{LeTaAA} for the precise definition of ``rank."
Roughly speaking, the rank of a complex reflection group $W$
is the minimal $n$ such that
$W$ can be realized as reflection group on $\mathbb C^n$.}
Well-generation can be equivalently characterised by a duality
property due to Orlik and Solomon \cite{OS}. Namely, a
complex reflection group of rank $n$ has two sets of distinguished integers
$d_1\leq d_2\leq \cdots \leq d_n$ and $d_1^*\geq d_2^*\geq \cdots \geq
d_n^*$, called its {\it degrees} and {\it codegrees}, respectively
(see \cite[p.~51 and Def.~10.27]{LeTaAA}). Orlik and 
Solomon observed, using case-by-case checking, that an irreducible
complex reflection group $W$ of rank $n$ is
well-generated if and only if its degrees and codegrees satisfy        
\begin{equation*}
d_i+d_i^*=d_n
\end{equation*}
for all $i=1,2,\dots,n$. The reader is referred to 
\cite[App.~D.2]{LeTaAA} for a table of
the degrees and codegrees of all irreducible complex reflection
groups. Together with the classification of 
Shephard and Todd \cite{ShToAA}, this constitutes a classification of 
well-generated complex reflection groups: 
the irreducible
well-generated complex reflection groups are 
\begin{enumerate}
\item[---]the two infinite
families $G(d,1,n)$ and $G(e,e,n)$, where $d,e,n$ are positive
integers, 
\item[---]the exceptional groups 
$G_4,G_5,G_6,G_8,G_9,G_{10},G_{14},G_{16},G_{17},G_{18},G_{20},G_{21}$ 
of\break rank~$2$,
\item[---]the exceptional groups 
$G_{23}=H_3,G_{24},G_{25},G_{26},G_{27}$ of rank $3$,
\item[---]the exceptional groups 
$G_{28}=F_4,G_{29},G_{30}=H_4,G_{32}$ of rank $4$,
\item[---]the exceptional group
$G_{33}$ of rank $5$,
\item[---]the exceptional groups 
$G_{34},G_{35}=E_6$ of rank $6$,
\item[---]the exceptional group 
$G_{36}=E_7$ of rank $7$,
\item[---]and the exceptional group 
$G_{37}=E_8$ of rank $8$.
\end{enumerate}

\noindent
In this list, we have made visible the groups
$H_3,F_4,H_4,E_6,E_7,E_8$ which appear as exceptional groups
in the classification of all
irreducible {\it real\/} reflection groups (cf.\ \cite{HumpAC}).

Let $W$ be a well-generated complex reflection group of rank $n$,
and let $T\subseteq W$ denote the set of {\it all\/} 
(complex) reflections in the group. 
Let $\ell_T:W\to\mathbb{Z}$ denote the word length in terms of the 
generators $T$. This word length is called {\it absolute length\/} or
{\it reflection length}. Furthermore,
we define a partial order $\le_T$ on $W$ by
\begin{equation} \label{eq:absord} 
u\le_T w\quad \text{if and only if}\quad
\ell_T(w)=\ell_T(u)+\ell_T(u^{-1}w).
\end{equation}
This partial order is called {\it absolute order} or {\it reflection
order}. As is well-known and easy to see, 
the equation in \eqref{eq:absord} is equivalent to the statement that every
shortest representation of $u$ by reflections
occurs as an initial segment in some shortest product representation
of $w$ by reflections. 

Now fix a (generalised) Coxeter 
element\footnote{An element of an irreducible well-generated complex
reflection group $W$ of rank $n$ is called a
{\it Coxeter element} if it is {\it regular} in the sense of
Springer \cite{SpriAA} (see also \cite[Def.~11.21]{LeTaAA}) 
and of order $d_n$. An element of $W$ is called regular if it
has an eigenvector which lies in no reflecting
hyperplane of a reflection of $W$. It follows from 
an observation of Lehrer and Springer, proved uniformly by Lehrer and 
Michel \cite{LeMiAA} (see \cite[Theorem~11.28]{LeTaAA}), that there is
always a regular element of order $d_n$
in an irreducible well-generated complex reflection group $W$ of rank $n$.
More generally, if a well-generated complex reflection group $W$
decomposes as $W\cong W_1\times W_2\times\dots\times W_k$, where the
$W_i$'s are irreducible, then a Coxeter element of $W$ is an element
of the form $c=c_1c_2\cdots c_k$, where $c_i$ is a Coxeter element of
$W_i$, $i=1,2,\dots,k$.
If $W$ is a {\it
real\/} reflection group, that is, if all generators in $T$ have order
$2$, then the notion of generalised Coxeter element given above 
reduces to that of a Coxeter element in the classical sense
(cf.\ \cite[Sec.~3.16]{HumpAC}).} $c\in W$ 
and a positive integer $m$. 
The {\it $m$-divisible non-crossing 
partitions} $NC^m(W)$ are defined as the set
\begin{multline*}
NC^m(W)=\big\{(w_0;w_1,\dots,w_m):w_0w_1\cdots w_m=c\text{ and }\\
\ell_T(w_0)+\ell_T(w_1)+\dots+\ell_T(w_m)=\ell_T(c)\big\}.
\end{multline*}
A partial order is defined on this set by 
$$
(w_0;w_1,\ldots,w_m) \leq (u_0;u_1,\ldots,u_m) \quad \text{if and only if}\quad u_i\le_T w_i \text{ for } 1\leq i\leq m.
$$
We have suppressed the dependence on $c$, since we understand this
definition up to isomorphism of posets. To be more precise,
it can be shown that any two Coxeter elements are related to each other
by conjugation and (possibly) an automorphism on the field of complex
numbers
(see \cite[Theorem~4.2]{SpriAA} or \cite[Cor.~11.25]{LeTaAA}), 
and hence the resulting posets $NC^m(W)$ are isomorphic to each other. 
If $m=1$, then $NC^1(W)$ can be identified with the set
$NC(W)$ of 
non-crossing partitions for the (complex) reflection group $W$ as defined
by Bessis and Corran (cf.\ \cite{BeCoAA} and
\cite[Sec.~13]{BesDAB}; their definition 
extends the earlier definition by Bessis
\cite{BesDAA} and Brady and Watt \cite{BRWaAA} for real reflection
groups).

The following result has been proved by a collaborative
effort of several authors (see \cite[Prop.~13.1]{BesDAB}).

\begin{theorem} \label{thm:2}
Let $W$ be an irreducible well-generated complex 
reflection group, and let
$d_1\le d_2\le\dots\le d_n$ be its degrees and $h:=d_n$ its Coxeter number. 
Then
\begin{equation} \label{eq:F-C}
\vert NC^m(W)\vert=\prod_{i=1}^n \frac {mh+d_i} {d_i}.
\end{equation}
\end{theorem}

\begin{remark} \label{rem:0}
(1)
The number in \eqref{eq:F-C} is called the {\it Fu\ss--Catalan number}
for the reflection group $W$.

\smallskip
(2) If $c$ is a Coxeter element of a well-generated complex
reflection group $W$ of rank $n$, then $\ell_T(c)=n$.
(This follows from \cite[Sec.~7]{BesDAB}.)
\end{remark}

We conclude this section by recalling the definition of decomposition
numbers from \cite{KratCB,KratCF,KrMuAB}. Although we need them here
only for (very small) real reflection groups, and although, strictly
speaking, they have been only defined for real reflection groups in
\cite{KratCB,KratCF,KrMuAB}, this definition can be extended to
well-generated complex reflection groups without any extra effort, 
which we do now.

Given a well-generated complex reflection group $W$ of rank $n$, types
$T_1,T_2,\dots,T_d$ (in the sense of the classification of well-generated 
complex reflection groups) such that the sum of the ranks of the
$T_i$'s equals $n$, and a Coxeter element $c$,
the {\it decomposition number} $N_W(T_1,T_2,\dots,T_d)$ is defined
as the number of ``minimal" factorisations $c=c_1c_2\cdots c_d$,
``minimal" meaning that 
$\ell_T(c_1)+\ell_T(c_2)+\dots+\ell_T(c_d)=\ell_T(c)=n$, 
such that, for $i=1,2,\dots,d$, the
type of $c_i$ as a parabolic Coxeter element is $T_i$.
(Here, the term ``parabolic Coxeter element" means a Coxeter element
in some parabolic subgroup. It follows from \cite[Prop.6.3]{RipoAA}
that any element $c_i$ is indeed a Coxeter element in a 
unique parabolic subgroup
of $W$.\footnote{The uniqueness can be argued as follows: suppose that
$c_i$ were a Coxeter element in two parabolic subgroups of $W$,
say $U_1$ and $U_2$. Then it must also be a Coxeter element in the
intersection $U_1\cap U_2$. On the other hand, the absolute length
of a Coxeter element of a complex reflection group $U$ 
is always equal to $\rk(U)$, the rank of $U$.
(This follows from the fact that, for each element
$u$ of $U$, we have
$\ell_T(u)=\text{codim}\big(\text{ker}(u-\text{id})\big)$, with id denoting
the identity element in $U$; see e.g.\ \cite[Prop.~1.3]{RipoAA}).
We conclude that $\ell_T(c_i)=\rk(U_1)=\rk(U_2)=\rk(U_1\cap U_2)$,
This implies that $U_1=U_2$.} 
By definition, the type of $c_i$ is the type of this
parabolic subgroup.) Since any two Coxeter elements are related 
to each other by
conjugation plus field automorphism, the decomposition numbers are independent
of the choice of the Coxeter element $c$.

The decomposition numbers for real reflection groups have been
computed in \cite{KratCB,KratCF,KrMuAB}. To compute the decomposition
numbers for well-generated
complex reflection groups is a task that remains to be done.

\section{Cyclic sieving I}
\label{sec:siev1}

In this section we present the first cyclic sieving conjecture due to
Armstrong \cite[Conj.~5.4.7]{ArmDAA}, and to Bessis and Reiner
\cite[Conj.~6.4]{BeReAA}.

Let $\phi:NC^m(W)\to NC^m(W)$ be the map defined by
\begin{equation} \label{eq:phi}
(w_0;w_1,\dots,w_m)\mapsto
\big((cw_mc^{-1})w_0(cw_mc^{-1})^{-1};
cw_mc^{-1},w_1,w_2,\dots,w_{m-1}\big).
\end{equation}
It is indeed not difficult to see that, if the $(m+1)$-tuple
on the left-hand side is an element of $NC^m(W)$, then so is
the $(m+1)$-tuple on the right-hand side.
For $m=1$, this action reduces to conjugation by the Coxeter element
$c$ (applied to $w_1$). Cyclic sieving arising from conjugation by $c$ 
has been the subject of \cite{BeReAA}.

It is easy to see that $\phi^{mh}$ acts as the identity,
where $h$ is the Coxeter number of $W$ (see \eqref{eq:Aktion}
and Lemma~\ref{lem:4} below).
By slight abuse of notation, let $C_1$ be the cyclic group of order $mh$
generated by $\phi$. (The slight abuse consists in the fact that we
insist on $C_1$ to be a cyclic group of order $mh$, while it may
happen that the order of the action of $\phi$ given in \eqref{eq:phi}
is actually a proper divisor of $mh$.)

Given these definitions, we are now in the position to state
the first cyclic sieving conjecture of
Armstrong, respectively of Bessis and Reiner.
By the results of \cite{KratCG} and of this paper, it becomes the
following theorem.

\begin{theorem} \label{thm:1}
For an irreducible well-generated complex reflection group $W$ and any $m\ge1$, 
the triple $(NC^m(W),\Cat^m(W;q),C_1)$, where $\Cat^m(W;q)$ is
the $q$-analogue of the Fu\ss--Catalan number defined by 
\begin{equation} \label{eq:FCZahl}
\Cat^m(W;q):=\prod_{i=1}^n \frac {[mh+d_i]_q} {[d_i]_q},
\end{equation}
exhibits the cyclic sieving phenomenon in the sense of
Reiner, Stanton and White \cite{ReSWAA}.
Here, $n$ is the rank of $W$, $d_1,d_2,\dots,d_n$ are the
degrees of $W$, $h$ is the Coxeter number of $W$, 
and $[\alpha]_q:=(1-q^\alpha)/(1-q)$.
\end{theorem}

\begin{remark}
We write $\Cat^m(W)$ for $\Cat^m(W;1)$.
\end{remark}

By definition of the cyclic sieving phenomenon, we have to
prove that $\Cat^m(W;q)$ is a polynomial in $q$ with non-negative
integer coefficients, and that
\begin{equation} \label{eq:1}
\vert\Fix_{NC^m(W)}(\phi^{p})\vert = 
\Cat^m(W;q)\big\vert_{q=e^{2\pi i p/mh}}, 
\end{equation}
for all $p$ in the range $0\le p<mh$.
The first fact is established in the next section, while the 
proof of the second is achieved by making use of 
several auxiliary results, given
in Section~\ref{sec:aux1}, to reduce the proof to a finite problem,
and a subsequent case-by-case analysis.
All details of this analysis can be found in \cite[Sec.~6]{KrMuAE}.
In the present paper, 
we content ourselves with discussing the cases where $W=G_{24}$
and where $W=G_{37}=E_8$, since these suffice to convey the flavour of
the necessary computations.

\section{The $q$-Fusz--Catalan numbers $\Cat^m(W;q)$}
\label{sec:pol}

The purpose of this section is to provide an elementary, self-contained
proof of the fact that, for all irreducible complex reflection groups
$W$, the $q$-Fu\ss--Catalan number $\Cat^m(W;q)$
is a polynomial in $q$ with non-negative integer coefficients.
For most of the groups, this is a known property. However, aside from the
fact that, for many of the known cases, the proof is very indirect and uses
deep algebraic results on rational Cherednik algebras, there still remained
some cases where this property had not been formally established.
The reader is referred to the ``Theorem" in Section~1.6 of \cite{GoGrAA},
which says that, under the assumption of a certain rank condition
(\cite[Hypothesis~2.4]{GoGrAA}), the $q$-Fu\ss--Catalan number
$\Cat^m(W;q)$ is a Hilbert series of a 
finite-dimensional quotient of the ring of invariants
of $W$ and also the graded character of a finite-dimensional irreducible
representation of a spherical rational Cherednik algebra associated with
$W$. At present, this rank condition has been proven
for all irreducible well-generated complex reflection groups apart
from $G_{17},G_{18},G_{29},G_{33},G_{34}$; see \cite[Tables~8 and 9, 
column~``rank"]{MaMiAA}, and the recent paper \cite{MariAA},
which establishes the result in the case of $G_{32}$.

In the sequel, aside from the standard notation
$[\al]_q=(1-q^\al)/(1-q)$ for $q$-integers, we shall also use the
$q$-binomial coefficient, which is defined by
$$
\begin{bmatrix} n\\k\end{bmatrix}_q:=\begin{cases} 1,&\text{if $k=0$,}\\
\frac {[n]_q\,[n-1]_q\cdots [n-k+1]_q} {[k]_q\,[k-1]_q\cdots [1]_q},
&\text{if $k>0$.}\end{cases}
$$

We begin with several auxiliary results.

\begin{proposition} \label{prop:1}
For all non-negative integers $n$ and $k$,
the $q$-binomial coefficient $\left[\begin{smallmatrix}
n\\k \end{smallmatrix}\right]_q$
is a polynomial in $q$ with non-negative
integer coefficients.
\end{proposition}

\begin{proof}This is a well-known fact, which can be derived either
from the recurrence relation(s) satisfied by the $q$-binomial
coefficients  (generalising Pascal's recurrence
relation for binomial coefficients;
cf.\ \cite[eqs.~(3.3.3) and (3.3.4)]{AndrAF}), 
or from the fact that the
$q$-binomial coefficient $\left[\begin{smallmatrix}
n\\k \end{smallmatrix}\right]_q$ is the generating function for
(integer) partitions with at most $k$ parts all of which are at
most $n-k$ (cf.\ \cite[Theorem~3.1]{AndrAF}).
\end{proof}

\begin{proposition} \label{prop:2}
For all non-negative integers $m$ and $n$, the $q$-Fu\ss--Catalan
number of type $A_n$,
$$
\frac {1} {[(m+1)n+1]_q}\begin{bmatrix} (m+1)n+1\\n\end{bmatrix}_q,
$$
is a polynomial in $q$ with non-negative
integer coefficients.
\end{proposition}

\begin{proof}
In \cite[Sec.~3.3]{LoehAA}, Loehr proves that
\begin{multline} \label{eq:Loehr}
\frac {1} {[(m+1)n+1]_q}\begin{bmatrix} (m+1)n+1\\n\end{bmatrix}_q\\
=\sum_{v\in\mathcal V_n^{(m)}}
q^{m\binom n2+\sum_{i\ge0}\left(m\binom {v_i}2-iv_i\right)}
\prod_{i\ge1}q^{v_i\sum_{j=1}^m (m-j)v_{i-j}}
\begin{bmatrix} v_i+v_{i-1}+\dots+v_{i-m}-1\\v_i\end{bmatrix}_q,
\end{multline}
where $\mathcal V_n^{(m)}$ denotes the set of all sequences
$v=(v_0,v_1,\dots,v_s)$ (for some $s$) of non-negative integers with
$v_0>0$, $v_s>0$, and $v_0+v_1+\dots+v_s=n$, and such that there is
never a string of $m$ or more
consecutive zeroes in $v$. By convention, $v_i = 0$ for all negative
$i$. His proof works by showing that the expressions on both sides of
\eqref{eq:Loehr} satisfy the same recurrence relation and initial
conditions, using classical $q$-binomial identities. We refer the
reader to \cite{LoehAA} for details. By Proposition~\ref{prop:1},
the expression on the right-hand side of \eqref{eq:Loehr} is
manifestly a polynomial in $q$ with non-negative integer coefficients.
\end{proof}

\begin{lemma} \label{lem:A}
If $a$ and $b$ are coprime positive integers, then
\begin{equation} \label{eq:ab} 
\frac {[ab]_q} {[a]_q\,[b]_q}
\end{equation}
is a polynomial in $q$ of degree $(a-1)(b-1)$, 
all of whose coefficients are in $\{0,1,-1\}$.
Moreover, if one disregards the coefficients which are $0$,
then $+1$'s and $(-1)$'s alternate, and the constant coefficient
as well as the leading coefficient of the polynomial equal $+1$.
\end{lemma}

\begin{proof}
Let $\Phi_n(q)$ denote the $n$-th cyclotomic polynomial in $q$.
Using the classical formula
$$1-q^n=
\prod _{d\mid n} ^{}\Phi_d(q),$$
we see that
$$
\frac {(1-q)(1-q^{ab})} {(1-q^a)(1-q^b)}=
\underset{{d_2\mid a,\,d_2\ne1}}{\prod _{d_1\mid a,\,d_1\ne1} ^{}}
\Phi_{d_1d_2}(q),
$$
so that, manifestly, the expression in \eqref{eq:ab} is a polynomial
in $q$. The claim concerning the degree of this polynomial is obvious.

In order to establish the claim on the coefficients,
we start with a sub-expression of \eqref{eq:ab},
\begin{equation} \label{eq:ab0}
\frac {(1-q^{ab})} {(1-q^a)(1-q^b)}=
\bigg(\sum_{i=0}^{b-1}q^{ia}\bigg)
\bigg(\sum_{j=0}^{\infty}q^{jb}\bigg)=
\sum_{k=0}^{\infty}C_kq^k,
\end{equation}
say. The assumption that $a$ and $b$ are coprime
implies that $0\le C_k\le1$ for  $k\le (a-1)(b-1)$.
Multiplying both sides of \eqref{eq:ab0} by $1-q$, we obtain
the equation
\begin{equation} \label{eq:ab1} 
\frac {[ab]_q} {[a]_q\,[b]_q}
=(1-q)\sum_{k=0}^{(a-1)(b-1)}C_kq^k
+(1-q)\sum_{k=(a-1)(b-1)+1}^{\infty}C_kq^k.
\end{equation}
By our previous observation on the coefficients $C_k$ with
$k\le(a-1)(b-1)$, it is obvious that the coefficients of the first
expression on the right-hand side of \eqref{eq:ab1} are alternately
$+1$ and $-1$, when $0$'s are disregarded. Since we already know
that the left-hand side is a polynomial in $q$ of degree $(a-1)(b-1)$, 
we may ignore the second expression.

The proof is concluded by observing that the claims on the constant
and leading coefficients are obvious.
\end{proof}

\begin{corollary} \label{cor:A}
Let $a$ and $b$ be coprime positive integers, and let $\ga$ be an integer
with $\ga\ge (a-1)(b-1)$. Then the expression
$$
\frac {[\ga]_q\,[ab]_q} {[a]_q\,[b]_q}
$$
is a polynomial in $q$ with non-negative
integer coefficients.
\end{corollary}

\begin{proof}Let
$$
\frac {[ab]_q} {[a]_q\,[b]_q}
=\sum_{k=0}^{(a-1)(b-1)}D_kq^k.
$$
We then have
\begin{equation} \label{eq:abc} 
\frac {[\ga]_q\,[ab]_q} {[a]_q\,[b]_q}
=\sum_{N=0}^{(a-1)(b-1)+\ga-1}q^N\sum_{k=\max\{0,N-\ga+1\}}^{N}D_k.
\end{equation}
If $N\le \ga-1$, then, by Lemma~\ref{lem:A}, 
the sum over $k$ on the right-hand side of \eqref{eq:abc} equals
$1-1+1-1+\cdots$, which is manifestly non-negative. 
On the other hand, if $N> \ga-1$, then we may rewrite
the sum over $k$ on the right-hand side of \eqref{eq:abc} as
$$
\sum_{k=\max\{0,N-\ga+1\}}^{N}D_k=
\sum_{k=N-\ga+1}^{(a-1)(b-1)}D_k=
\sum_{k=0}^{(a-1)(b-1)+\ga-1-N}D_{(a-1)(b-1)-k}.
$$
Again, by Lemma~\ref{lem:A}, 
this sum equals
$1-1+1-1+\cdots$, which is manifestly non-negative. 
\end{proof}

The next lemmas all have a very similar flavour, and so do their
proofs. In order to avoid repetition, proof details are only
provided for Lemmas~\ref{lem:B} and \ref{lem:61015}; the
proofs of Lemmas~\ref{lem:C}--\ref{lem:J}, \ref{lem:679}--\ref{lem:4710}
follow the pattern exhibited in the proof of Lemma~\ref{lem:B},
while the proofs of Lemmas~\ref{lem:61015B}--\ref{lem:101215C}
follow that of the proof of Lemma~\ref{lem:J}.
Full details are found in \cite[Sec.~4]{KrMuAE}.

\begin{lemma} \label{lem:B}
Let $\al$ and $\be$ be positive integers with $\al\ge 6$ and $\be\ge
8$. Then the expression
$$
\left[\al\right]_{q^3}\left[\be\right]_{q^4}
\frac {\left[72\right]_q\left[3\right]_q\left[4\right]_q} 
{\left[8\right]_q\left[9\right]_q\left[12\right]_q}
$$
is a polynomial in $q$ with non-negative
integer coefficients.
\end{lemma}

\begin{proof}
We have
\begin{multline*}
\frac {\left[72\right]_q\left[3\right]_q\left[4\right]_q} 
{\left[8\right]_q\left[9\right]_q\left[12\right]_q}\\
=(1
-q^3
+q^9
-q^{15}
+q^{18}
)
(1
-q^4
+q^8
-q^{12}
+q^
   {16}
-q^{20}
+q^{24}
-q^{28}
+q^{32}
).
\end{multline*}
It should be observed that both factors on the right-hand side
have the property that coefficients are in $\{0,1,-1\}$ and that
$(+1)$'s and $(-1)$'s alternate, if one disregards the coefficients
which are $0$.
If we now apply the same idea as in the proof of Corollary~\ref{cor:A},
then we see that $[\al]_{q^3}$ times the first factor 
is a polynomial in $q$ with
non-negative integer coefficients, as is
$[\be]_{q^4}$ times the second factor.
Taken together, this establishes
the claim.
\end{proof}

\begin{lemma} \label{lem:G}
Let $\al$ and $\be$ be positive integers with $\al\ge 26$ and $\be\ge
8$. Then the expression
$$
\left[\al\right]_{q}\left[\be\right]_{q^4}
\frac {\left[15\right]_q}
{\left[3\right]_q\left[5\right]_q}
\frac {\left[72\right]_q\left[3\right]_q\left[4\right]_q} 
{\left[8\right]_q\left[9\right]_q\left[12\right]_q}
$$
is a polynomial in $q$ with non-negative
integer coefficients.
\end{lemma}

\begin{lemma} \label{lem:C}
Let $\al$ and $\be$ be positive integers with $\al\ge 18$ and $\be\ge
3$. Then the expression
$$
\left[\al\right]_{q^3}\left[\be\right]_{q^4}
\frac {\left[90\right]_q\left[3\right]_q\left[4\right]_q} 
{\left[5\right]_q\left[6\right]_q\left[9\right]_q}
$$
is a polynomial in $q$ with non-negative
integer coefficients.
\end{lemma}

\begin{lemma} \label{lem:F}
Let $\al$ and $\be$ be positive integers with $\al\ge 20$ and $\be\ge
18$. Then the expression
$$
\left[\al\right]_{q}\left[\be\right]_{q^3}
\frac {\left[90\right]_q\left[3\right]_q}
{\left[5\right]_q\left[6\right]_q\left[9\right]_q}
$$
is a polynomial in $q$ with non-negative
integer coefficients.
\end{lemma}

\begin{lemma} \label{lem:D}
Let $\al$ be a positive integer with $\al\ge 26$.
Then the expression
$$
\left[\al\right]_{q}
\frac {\left[15\right]_q}
{\left[3\right]_q\left[5\right]_q}
\frac {\left[12\right]_{q^3}}
{\left[3\right]_{q^3}\left[4\right]_{q^3}}
$$
is a polynomial in $q$ with non-negative
integer coefficients.
\end{lemma}

\begin{lemma} \label{lem:E}
Let $\al$ be a positive integer with $\al\ge 14$.
Then the expression
$$
\left[\al\right]_{q}
\frac {\left[15\right]_q}
{\left[3\right]_q\left[5\right]_q}
\frac {\left[6\right]_{q^3}}
{\left[2\right]_{q^3}\left[3\right]_{q^3}}
$$
is a polynomial in $q$ with non-negative
integer coefficients.
\end{lemma}

\begin{lemma} \label{lem:H}
Let $\al$ and $\be$ be positive integers with $\al\ge 30$ and $\be\ge
20$. Then the expression
$$
\left[\al\right]_{q}\left[\be\right]_{q^2}
\frac {\left[84\right]_q\left[2\right]_q}
{\left[4\right]_q\left[6\right]_q\left[7\right]_q}
$$
is a polynomial in $q$ with non-negative
integer coefficients.
\end{lemma}

\begin{lemma} \label{lem:I}
Let $\al$ and $\be$ be positive integers with $\al\ge 24$ and $\be\ge
68$. Then the expression
$$
\left[\al\right]_{q}\left[\be\right]_{q}
\frac {\left[105\right]_q}
{\left[3\right]_q\left[5\right]_q\left[7\right]_q}
$$
is a polynomial in $q$ with non-negative
integer coefficients.
\end{lemma}

\begin{lemma} \label{lem:J}
Let $\al$ and $\be$ be positive integers with $\al\ge 24$ and $\be\ge
34$. Then the expression
$$
\left[\al\right]_{q}\left[\be\right]_{q}
\frac {\left[70\right]_q}
{\left[2\right]_q\left[5\right]_q\left[7\right]_q}
$$
is a polynomial in $q$ with non-negative
integer coefficients.
\end{lemma}

\begin{lemma} \label{lem:61015}
Let $\al$ and $\be$ be positive integers with $\al\ge 4$ and $\be\ge
2$. Then the expression
$$
\left[\al\right]_{q^2}\left[\be\right]_{q^5}
\frac {\left[30\right]_q\left[2\right]_q\left[3\right]_q\left[5\right]_q}
{\left[6\right]_q\left[10\right]_q\left[15\right]_q}
$$
is a polynomial in $q$ with non-negative
integer coefficients.
\end{lemma}

\begin{proof}
We have
$$
\frac {\left[30\right]_q\left[2\right]_q\left[3\right]_q\left[5\right]_q}
{\left[6\right]_q\left[10\right]_q\left[15\right]_q}
=
1 + q - q^3 - q^4 - q^5 + q^7 + q^8.
$$
If we multiply this expression by $[\al]_{q^2}$, then, for $\al=4$ we obtain
$$
1 + q + q^2 - q^5 - q^9 + q^{12} + q^{13} + q^{14},
$$
for $\al=5$ we obtain
$$
1 + q + q^2 - q^5 + q^8 - q^{11} + q^{14} + q^{15} + q^{16},
$$
and, for $\al\ge6$, we obtain 
$$
1 + q + q^2 - q^5 + 
q^8+q^{10}+p_1(q)+q^{2\al-4}+q^{2\al-2}
 - q^{2\al+1} + q^{2\al+4} + q^{2\al+5} + q^{2\al+6},
$$
where $p_1(q)$ is a polynomial
in $q$ with non-negative coefficients of order 
at least $11$ and degree at most $2\al-5$. 
In all cases it is obvious that the product of the result
and $[\be]_{q^5}$, with $\be\ge2$, 
is a polynomial in $q$ with non-negative coefficients.
\end{proof}

\begin{lemma} \label{lem:61015B}
Let $\al$ and $\be$ be positive integers with $\al\ge 14$ and $\be\ge
2$. Then the expression
$$
\left[\al\right]_{q}\left[\be\right]_{q^5}
\frac {\left[14\right]_q}
{\left[2\right]_q\left[7\right]_q}
\frac {\left[30\right]_q\left[2\right]_q\left[3\right]_q\left[5\right]_q}
{\left[6\right]_q\left[10\right]_q\left[15\right]_q}
$$
is a polynomial in $q$ with non-negative
integer coefficients.
\end{lemma}

\begin{lemma} \label{lem:61015C}
Let $\al$ and $\be$ be positive integers with $\al\ge 32$ and $\be\ge
12$. Then the expression
$$
\left[\al\right]_{q}\left[\be\right]_{q^2}
\frac {\left[35\right]_q}
{\left[5\right]_q\left[7\right]_q}
\frac {\left[30\right]_q\left[2\right]_q\left[3\right]_q\left[5\right]_q}
{\left[6\right]_q\left[10\right]_q\left[15\right]_q}
$$
is a polynomial in $q$ with non-negative
integer coefficients.
\end{lemma}

\begin{lemma} \label{lem:101215}
Let $\al$ and $\be$ be positive integers with $\al\ge 16$ and $\be\ge
2$. Then the expression
$$
\left[\al\right]_{q^2}\left[\be\right]_{q^5}
\frac {\left[60\right]_q\left[2\right]_q\left[3\right]_q\left[5\right]_q}
{\left[10\right]_q\left[12\right]_q\left[15\right]_q}
$$
is a polynomial in $q$ with non-negative
integer coefficients.
\end{lemma}

\begin{lemma} \label{lem:101215B}
Let $\al$ and $\be$ be positive integers with $\al\ge 56$ and $\be\ge
4$. Then the expression
$$
\left[\al\right]_{q}\left[\be\right]_{q^2}
\frac {\left[35\right]_q}
{\left[5\right]_q\left[7\right]_q}
\frac {\left[60\right]_q\left[2\right]_q\left[3\right]_q\left[5\right]_q}
{\left[10\right]_q\left[12\right]_q\left[15\right]_q}
$$
is a polynomial in $q$ with non-negative
integer coefficients.
\end{lemma}

\begin{lemma} \label{lem:101215C}
Let $\al$ and $\be$ be positive integers with $\al\ge 38$ and $\be\ge
2$. Then the expression
$$
\left[\al\right]_{q}\left[\be\right]_{q^5}
\frac {\left[14\right]_q}
{\left[2\right]_q\left[7\right]_q}
\frac {\left[60\right]_q\left[2\right]_q\left[3\right]_q\left[5\right]_q}
{\left[10\right]_q\left[12\right]_q\left[15\right]_q}
$$
is a polynomial in $q$ with non-negative
integer coefficients.
\end{lemma}

\begin{lemma} \label{lem:679}
Let $\al$ and $\be$ be positive integers with $\al\ge 30$ and $\be\ge
26$. Then the expression
$$
\left[\al\right]_{q}\left[\be\right]_{q^3}
\frac {\left[126\right]_q\left[3\right]_q}
{\left[6\right]_q\left[7\right]_q\left[9\right]_q}
$$
is a polynomial in $q$ with non-negative
integer coefficients.
\end{lemma}

\begin{lemma} \label{lem:7912}
Let $\al$ and $\be$ be positive integers with $\al\ge 66$ and $\be\ge
54$. Then the expression
$$
\left[\al\right]_{q}\left[\be\right]_{q^3}
\frac {\left[252\right]_q\left[3\right]_q}
{\left[7\right]_q\left[9\right]_q\left[12\right]_q}
$$
is a polynomial in $q$ with non-negative
integer coefficients.
\end{lemma}

\begin{lemma} \label{lem:4710}
Let $\al$ and $\be$ be positive integers with $\al\ge 54$ and $\be\ge
34$. Then the expression
$$
\left[\al\right]_{q}\left[\be\right]_{q^2}
\frac {\left[140\right]_q\left[2\right]_q}
{\left[4\right]_q\left[7\right]_q\left[10\right]_q}
$$
is a polynomial in $q$ with non-negative
integer coefficients.
\end{lemma}

We are now ready for the proof of the main result of this section.

\begin{theorem} \label{thm:0}
For all irreducible well-generated complex reflection groups
and positive integers $m$,
the $q$-Fu\ss--Catalan number $\Cat^m(W;q)$ 
is a polynomial in $q$ with non-negative
integer coefficients.
\end{theorem}

\begin{proof}
First, let $W=A_n$. In this case, the degrees are $2,3,\dots,n+1$,
and hence
$$
\Cat^m(A_n;q)=\frac {1} {[(m+1)n+1]_q}\begin{bmatrix}
  (m+1)n+1\\n\end{bmatrix}_q, 
$$
which, by Proposition~\ref{prop:2}, is a polynomial in $q$ with 
non-negative integer coefficients.

Next, let $W=G(d,1,n)$. In this case, the degrees are 
$d,2d,\dots,nd$, and hence
$$
\Cat^m(G(d,1,n);q)=\begin{bmatrix}
  (m+1)n\\n\end{bmatrix}_{q^d}, 
$$
which, by Proposition~\ref{prop:1}, is a polynomial in $q$ with 
non-negative integer coefficients.

Now, let $W=G(e,e,n)$. In this case, the degrees are 
$e,2e,\dots,(n-1)e,n$, and hence
\begin{align*}
\Cat^m(G(e,e,n);q)&=
\frac {[m(n-1)e+n]_q} {[n]_q}
\prod _{i=1} ^{n-1}\frac {[m(n-1)e+ie]_q} {[ie]_q}\\
&=
\begin{bmatrix}
  (m+1)(n-1)\\n-1\end{bmatrix}_{q^e}+
q^n[e]_{q^n}\begin{bmatrix}
  (m+1)(n-1)\\n\end{bmatrix}_{q^e}, 
\end{align*}
which, by Proposition~\ref{prop:1}, is a polynomial in $q$ with 
non-negative integer coefficients.

\medskip
It remains to verify the claim for the exceptional groups.

\medskip
For the groups $W=G_6,G_9,G_{14},G_{17},G_{21},$
and partially for the groups $W=G_{20},G_{23},\break G_{28},
G_{30},G_{33},G_{35},G_{36},G_{37}$
(depending on congruence properties of the parameter $m$),
polynomiality and non-negativity of coefficients
of the corresponding $q$-Fu\ss--Catalan number can be directly
read off by a proper rearrangement of the terms in the defining
expression; for example, for $W=G_{21}$
(with degrees given by $12,60$) we have
$$
\Cat^m(G_{21};q)=\frac {[60m+12]_q\,[60m+60]_q} {[12]_q\,[60]_q}
=[5m+1]_{q^{12}}\,[m+1]_{q^{60}},
$$
which is manifestly a polynomial in $q$ with 
non-negative integer coefficients.

\medskip
For the groups $G_5,G_{10},G_{18},G_{26},G_{27},G_{29},G_{34}$,
the terms in the defining expression of the corresponding $q$-Fu\ss--Catalan
number can be arranged in a manner so that a $q$-binomial coefficient
appears; polynomiality and non-negativity of coefficients then follow
from Proposition~\ref{prop:1}. For example, for $W=G_{34}$ 
(with degrees given by $6,12,18,24,30,42$) we have
\begin{align*}
\Cat^m(G_{34};q)&=\frac {[42m+6]_q\,[42m+12]_q\,[42m+18]_q\,
[42m+24]_q\,[42m+30]_q\,[42m+42]_q} 
{[6]_q\,[12]_q\,[18]_q\,[24]_q\,[30]_q\,[42]_q}\\
&=[m+1]_{q^{42}}
\begin{bmatrix} 7m+5\\5\end{bmatrix}_{q^{6}},
\end{align*}
which, written in this form, is
obviously a polynomial in $q$ with non-negative integer
coefficients.

\medskip
On the other hand,
for the groups $G_4,G_{8},G_{16},G_{25},G_{32}$,
the terms in the defining expression of the corresponding $q$-Fu\ss--Catalan
number can be arranged in a manner so that a $q$-Fu\ss--Catalan
number of type $A$ appears and Proposition~\ref{prop:2} applies;
for example, for $W=G_{32}$
(with degrees given by $12,18,24,30$) we have
\begin{align*}
\Cat^m(G_{32};q)&=\frac
    {[30m+12]_q\,[30m+18]_q\,[30m+24]_q\,[30m+30]_q} 
{[12]_q\,[18]_q\,[24]_q\,[30]_q}\\
&=
\frac {1} {[5m+6]_{q^{6}}}
\begin{bmatrix}
5m+6\\5\end{bmatrix}_{q^{6}}, 
\end{align*}
which indeed fits into the framework of Proposition~\ref{prop:2} and, hence, 
is a polynomial in $q$ with 
non-negative integer coefficients.

\medskip
In the other cases, the more ``specialised" auxiliary results
given in Corollary~\ref{cor:A} and 
Lemmas~\ref{lem:B}--\ref{lem:4710} have to be applied.
For the sake of illustration,
we exhibit one example for each of them below, with full details
being provided in \cite[Sec.~4]{KrMuAE}.
In general, the idea is that, given a rational expression
consisting of cyclotomic factors, as in the definition of the
$q$-Fu\ss--Catalan numbers, one tries to place denominator factors
below appropriate numerator factors so that one can divide out
the denominator factor completely. For example, if we were to
encounter the expression 
$$
\frac {[30m+12]_q\cdot \text{(other terms)}} 
{[12]_q\cdot \text{(other terms)}}
$$
and know that $m$ is even, then we would try to simplify 
this to
$$
\left[\tfrac {5m+2} {2}\right]_{q^{12}}\cdot
\frac {\text{(other terms)}} {\text{(other terms)}},
$$
where $[\tfrac {5m+2} {2}]_{q^{12}}$ is manifestly a polynomial in $q$
with non-negative integer coefficients.
On the other hand, in a situation where {\it two} denominator factors
``want" to divide a {\it single} numerator factor, we ``extract" as
much as we can from the numerator factor and compensate by 
additional ``fudge" factors. To be more concrete, if we
encounter the expression
$$
\frac {[14m+14]_q\cdot \text{(other terms)}} 
{[6]_q\left[14\right]_q\cdot \text{(other terms)}}
$$
and we know that $m\equiv0~(\text{mod }3)$, 
then we would try the rewriting
$$
\left[\tfrac {m+1} {3}\right]_{q^{42}}
\frac {[21]_{q^2}} {[3]_{q^2}\left[7\right]_{q^2}\left[2\right]_q}
\cdot
\frac {\text{(other terms)}} {\text{(other terms)}},
$$
with the idea that we might find somewhere else a term $[2\al]_{q}$,
which could be combined with the term $[2]_q$ in the denominator
into $[2\al]_q/[2]_q=[\al]_{q^2}$, and then apply
Corollary~\ref{cor:A} to see that 
$$
[\al]_{q^2}\frac {[21]_{q^2}} {[3]_{q^2}\left[7\right]_{q^2}}
$$
is a polynomial in $q$ with non-negative integer coefficients
(provided $\al$ is at least $12$), with 
$\left[\tfrac {m+1} {3}\right]_{q^{42}}$ being such a polynomial
in any case.

In situations where {\it three} denominator factors ``want" to divide
a {\it single} numerator factor, one has to perform more complicated 
rearrangements, in order to be able to apply one of the
Lemmas~\ref{lem:B}--\ref{lem:4710}.

\medskip
For example, 
for $W=G_{24}$, the degrees are $4,6,14$, and hence
$$
\Cat^m(G_{24};q)=\frac {[14m+4]_q\,[14m+6]_q\,[14m+14]_q} 
{[4]_q\,[6]_q\,[14]_q}.
$$
We have
$$
\Cat^m(G_{24};q)=
\begin{cases} 
\left[\tfrac {7m} {2}+1\right]_{q^4}\,
\left[\tfrac {14m} {6}+1\right]_{q^6}\,
\left[m+1\right]_{q^{14}},&\text{if $m\equiv0$ (mod 6),}\\
\left[\tfrac {7m+2} {3}\right]_{q^6}\,
\left[\tfrac {7m+3} {2}\right]_{q^4}\,
\left[m+1\right]_{q^{14}},&\text{if $m\equiv1$ (mod 6),}\\
\left[\tfrac {7m} {2}+1\right]_{q^4}\,
\left[7m+3\right]_{q^2}\,
\left[\frac {m+1} {3}\right]_{q^{42}}
\frac {[21]_{q^2}} {[3]_{q^2}\left[7\right]_{q^2}}
,&\text{if $m\equiv2$ (mod 6),}\\
\left[{7m}+ {2}\right]_{q^2}\,
\left[\frac{7m}3+1\right]_{q^6}\,
\left[\frac {m+1} {2}\right]_{q^{28}}
\frac {[14]_{q^2}} {[2]_{q^2}\left[7\right]_{q^2}}
,&\text{if $m\equiv3$ (mod 6),}\\
\left[\tfrac {7m+2} {6}\right]_{q^{12}}
\frac {[6]_{q^2}} {[2]_{q^2}\left[3\right]_{q^2}}
\left[7m+3\right]_{q^2}\,
\left[{m+1} \right]_{q^{14}}
,&\text{if $m\equiv4$ (mod 6),}\\
\left[ {7m} +{2}\right]_{q^2}\,
\left[\frac{7m+3}2\right]_{q^4}\,
\left[\frac {m+1} {3}\right]_{q^{42}}
\frac {[21]_{q^2}} {[3]_{q^2}\left[7\right]_{q^2}}
,&\text{if $m\equiv5$ (mod 6),}
\end{cases}
$$
which, by Corollary~\ref{cor:A}, are polynomials in $q$ with 
non-negative integer coefficients in all cases.

\medskip
For $W=G_{30}=H_4$, the degrees are $2,12,20,30$, and hence
$$
\Cat^m(H_4;q)=\frac {[30m+2]_q\,[30m+12]_q\,[30m+20]_q\,[30m+30]_q}
{[2]_q\,[12]_q\,[20]_q\,[30]_q}.
$$
If $m$ is odd, then we may write
\begin{align*}
\Cat^m(H_4;q)&=
\left[\tfrac{15m+1}2\right]_{q^{4}}\,
\left[5m+2\right]_{q^{6}}\,
\left[3m+2\right]_{q^{10}}\,
\left[\tfrac {m+1} {2}\right]_{q^{60}}\,
\frac {[30]_{q^2}\left[2\right]_{q^2}
\left[3\right]_{q^2}\left[5\right]_{q^2}} 
{[6]_{q^6}\left[10\right]_{q^2}\left[15\right]_{q^2}},
\end{align*}
which, by Lemma~\ref{lem:61015}, is a polynomial in $q$ with 
non-negative integer coefficients. 

\medskip
For $W=G_{35}=E_6$, the degrees are $2,5,6,8,9,12$, and hence
$$
\Cat^m(E_6;q)=\frac {[12m+2]_q\,[12m+5]_q\,[12m+6]_q\,[12m+8]_q\,
[12m+9]_q\,[12m+12]_q}
{[2]_q\,[5]_q\,[6]_q\,[8]_q\,[9]_q\,[12]_q}.
$$
If $m\equiv
5
~(\text{mod }30),$ then we have \begin{multline*}\Cat^m(E_6;q)=
{\left[ 6m+1 \right]_{q^{ 2}}} 
{\left[ \tfrac{12m+5}{5}  
   \right]_{q^{ 5}}} 
{\left[ 2 m+1 \right]_{q^{ 6}}} \\
\times
 {\left[ 3m+2 \right]_{q^{ 4 }}} {\left[ 
  {4m+3}  \right]_{q^{ 3 }}} 
{\left[\tfrac { m+1 }6
   \right]_{q^{ 72}}} 
\frac {\left[72\right]_q\left[3\right]_q\left[4\right]_q} 
{\left[8\right]_q\left[9\right]_q\left[12\right]_q}
,\end{multline*} 
which, by Lemma~\ref{lem:B}, 
is a polynomial in $q$ with 
non-negative integer coefficients.

If $m\equiv
7
~(\text{mod }30),$ then we have 
\begin{multline*}\Cat^m(E_6;q)=
{\left[\tfrac { 6m+1}2 \right]_{q^{ 4}}} 
{\left[ {12m+5}  
  \right]_{q}} {\left[ \tfrac {2 m+1}{15} \right]_{q^{ 90}}}\\
\times
\frac {\left[90\right]_q\left[3\right]_q\left[4\right]_q} 
{\left[5\right]_q\left[6\right]_q\left[9\right]_q}
 {\left[ 3m+2 \right]_{q^{ 4 }}} {\left[ 
  {4m+3} \right]_{q^{ 3 }}} {\left[ \tfrac {m+1}2 
   \right]_{q^{ 24}}} 
\frac {[6]_{q^4}} {[2]_{q^4}\left[3\right]_{q^4}}
,\end{multline*} 
which, by Corollary~\ref{cor:A} and Lemma~\ref{lem:C}, 
is a polynomial in $q$ with 
non-negative integer coefficients.

If $m\equiv
8
~(\text{mod }30),$ then we have 
\begin{multline*}\Cat^m(E_6;q)=
{\left[ 6m+1 \right]_{q^{ 2 }}} {\left[ 12m+5 
  \right]_{q}} {\left[ 2 m+1 \right]_{q^{ 6 
  }}} {\left[ \tfrac{3m+2}{2}  \right]_{q^{ 8 }}} \\
\times
 {\left[ \tfrac{4m+3}{5}  \right]_{q^{ 15 }}}
\frac {[15]_q} {[3]_q\left[5\right]_q} 
 {\left[ \tfrac{m+1}{3} \right]_{q^{ 36}}} \frac {{ 
  {\left[ 12 \right]_{q^{ 3}}} }} { {\left[ 3 
   \right]_{q^{ 3}}} {\left[ 4 \right]_{q^{ 3}}}}
,\end{multline*}
which, by Lemma~\ref{lem:D}, is a polynomial in $q$ with 
non-negative integer coefficients.

If $m\equiv
13
~(\text{mod }30),$ then we have \begin{multline*}\Cat^m(E_6;q)=
{\left[ 6m+1 \right]_{q^{ 2 }}} {\left[ 12m+5 
  \right]_{q}} {\left[ \tfrac{2m+1}{3}  
   \right]_{q^{ 18}}} \frac {{ {\left[ 6 \right]_{q^{ 3}}} 
  }} { {\left[ 2 \right]_{q^{ 3}}} {\left[ 3 
   \right]_{q^{ 3}}}} \\
\times{\left[ 3m+2 \right]_{q^4}} 
 {\left[ \tfrac{4m+3}{5}  \right]_{q^{ 15 }}} 
\frac {[15]_q} {[3]_q\left[5\right]_q}
 {\left[ \tfrac{m+1}{2} \right]_{q^{ 24}}} \frac {{ 
  {\left[ 6 \right]_{q^{ 4}}} }} { {\left[ 2 
   \right]_{q^{ 4}}} {\left[ 3 \right]_{q^{ 4}}}}
,\end{multline*} 
which, by Lemma~\ref{lem:E}, is a polynomial in $q$ with 
non-negative integer coefficients.

If $m\equiv
22
~(\text{mod }30),$ then we have \begin{multline*}\Cat^m(E_6;q)=
{\left[ 6m+1 \right]_{q^{ 2 }}} {\left[ 12m+5 
  \right]_{q}} {\left[ \tfrac{2 m+1}{15}  
   \right]_{q^{ 90}}} 
\frac {{ {\left[ 90 \right]_{q}} {\left[ 3 \right]_{q}} 
  }} { {\left[ 5 \right]_{q}} {\left[ 6 \right]_{q}} {\left[ 9 
   \right]_{q}}}\\
\times {\left[ \tfrac{3m+2}{2}  \right]_{q^{ 8 
  }}} {\left[ 4m+3 \right]_{q^3}} 
 {\left[ m+1 \right]_{q^{ 12}}} 
,\end{multline*} 
which, by Lemma~\ref{lem:F}, is a polynomial in $q$ with 
non-negative integer coefficients.

If $m\equiv
23
~(\text{mod }30),$ then we have \begin{multline*}\Cat^m(E_6;q)=
{\left[ 6m+1 \right]_{q^{ 2}}} 
{\left[ 12m+5
  \right]_{q}} {\left[ 2 m+1 \right]_{q^{ 6}}} \\
\times
 {\left[ {3m+2} \right]_{q^{ 4 }}} {\left[ 
   \tfrac{4m+3}5  \right]_{q^{ 15}}} 
\frac {[15]_q} {[3]_q\left[5\right]_q}
 {\left[ \tfrac {m+1}6 \right]_{q^{ 
   72}}} 
\frac {\left[72\right]_q\left[3\right]_q\left[4\right]_q} 
{\left[8\right]_q\left[9\right]_q\left[12\right]_q}
,\end{multline*} 
which, by Lemma~\ref{lem:G}, 
is a polynomial in $q$ with 
non-negative integer coefficients.

For $W=G_{36}=E_7$, the degrees are $2,6,8,10,12,14,18$, and hence
\begin{multline*}
\Cat^m(E_7;q)=\frac {[18m+2]_q\,[18m+6]_q\,[18m+8]_q\,[18m+10]_q\,
}
{[2]_q\,[6]_q\,[8]_q\,[10]_q}\\
\times
\frac {
[18m+12]_q\,[18m+14]_q\,[18m+18]_q}
{[12]_q\,[14]_q\,[18]_q}.
\end{multline*}
If $m\equiv
18
~(\text{mod }140),$ then we have \begin{multline*}\Cat^m(E_7;q)=
{\left[ 9 m+1 \right]_{q^{ 2}}} 
{\left[ \tfrac{3m+1}5 
   \right]_{q^{ 30}}} 
\frac {[15]_{q^2}} {[3]_{q^2}\left[5\right]_{q^2}}\\
\times
{\left[ \tfrac{9m+4}2 \right]_{q^{ 4}}} 
 {\left[ 9m+5 \right]_{q^{ 2}}} {\left[ \tfrac{3m+2}{28} 
   \right]_{q^{ 168}}} 
\frac {[84]_{q^2}\left[2\right]_{q^2}} 
{[4]_{q^2}\left[6\right]_{q^2}\left[7\right]_{q^2}}
{\left[ 9m+7 \right]_{q^{ 2}}} 
 {\left[ m+1 \right]_{q^{ 18}}}
,\end{multline*} 
which, by Corollary~\ref{cor:A} and Lemma~\ref{lem:H}, 
is a polynomial in $q$ with 
non-negative integer coefficients.

If $m\equiv
23
~(\text{mod }140),$ then we have \begin{multline*}\Cat^m(E_7;q)=
{\left[ \tfrac{9m+1}{4}  \right]_{q^{ 8}}} {\left[ 
   \tfrac{3 m+1}{35}  \right]_{q^{ 210}}} \frac {{ {\left[ 105
    \right]_{q^{ 2}}} }} { {\left[ 3 \right]_{q^{ 2}}}
{\left[ 5 \right]_{q^{ 2}}} 
  {\left[ 7 \right]_{q^{ 2}}} } {\left[ 9m+4 
  \right]_{q^{ 2}}} {\left[ 9 m+5 \right]_{q^2}}\\
\times {\left[ 
  3m+2 \right]_{q^{ 6}}} {\left[ 9 m+7 \right]_{q^2}} 
 {\left[ \tfrac{m+1}{2} \right]_{q^{ 36}}} \frac {{ 
  {\left[ 6 \right]_{q^{ 6}}} }} { {\left[ 2 
   \right]_{q^{ 6}}} {\left[ 3 \right]_{q^{ 6}}} }
,\end{multline*} 
which, by Corollary~\ref{cor:A} and Lemma~\ref{lem:I}, 
is a polynomial in $q$ with 
non-negative integer coefficients.

If $m\equiv
54
~(\text{mod }140),$ then we have \begin{multline*}\Cat^m(E_7;q)=
{\left[ 9 m+1 \right]_{q^{ 2}}} {\left[ 3 m+1 \right]_{q^{ 
  6}}} {\left[ \tfrac{9m+4}{70}  \right]_{q^{ 140}}} 
  \frac {{ {\left[ 70 \right]_{q^{ 2}}} }} { 
  {\left[ 2 \right]_{q^{ 2}}}
  {\left[ 5 \right]_{q^{ 2}}} {\left[ 7 \right]_{q^{ 2}}} 
  }
 {\left[ 9 m+5 \right]_{q^2}}\\
\times {\left[ 
   \tfrac{3m+2}{4}  \right]_{q^{ 24}}} \frac {{ {\left[ 6 
   \right]_{q^{ 4}}} }} { {\left[ 2 \right]_{q^{ 4}}} 
  {\left[ 3 \right]_{q^{ 4}}} } {\left[ 9 m+7 
  \right]_{q^2}} {\left[ m+1 \right]_{q^{ 18}}}
.\end{multline*} 
If one decomposes $[9m+7]_{q^2}$ as
$[\frac {9m} {2}+4]_{q^4}+q^2[\frac {9m} {2}+3]_{q^4}$,
then one sees that, by Corollary~\ref{cor:A} and Lemma~\ref{lem:J},
this is a polynomial in $q$ with 
non-negative integer coefficients.

For $W=G_{37}=E_8$, the degrees are $2,8,12,14,18,20,24,30$, and hence
\begin{multline*}
\Cat^m(E_7;q)=\frac {[30m+2]_q\,[30m+8]_q\,[30m+12]_q\,[30m+14]_q\,
}
{[2]_q\,[8]_q\,[12]_q\,[14]_q}\\
\times
\frac {
[30m+18]_q\,[30m+20]_q\,[30m+24]_q\,[30m+30]_q}
{[18]_q\,[20]_q\,[24]_q\,[30]_q}.
\end{multline*}
If $m\equiv
3
~(\text{mod }84),$ then we have 
\begin{multline*} \Cat^m(E_8;q)=
{\left[ \tfrac{15 m+1}2 \right]_{q^{ 4}}} 
{\left[ \tfrac{15m+4}7 \right]_{q^{ 
   14}}} 
 {\left[ 5m+2 \right]_{q^{ 6 }}}
 {\left[ \tfrac{15m+7}4 
   \right]_{q^{ 8}}} 
 {\left[ \tfrac{5m+3}6 \right]_{q^{ 36}}}
\frac {[6]_{q^6}} {[2]_{q^6}\left[3\right]_{q^6}}
\\
\times 
 {\left[ 3m+2 \right]_{q^{ 10 }}} {\left[ 5m+4
   \right]_{q^{ 6 }}} {\left[ \tfrac{m+1}4 \right]_{q^{ 120}}}
\frac {[60]_{q^2}\left[2\right]_{q^2}\left[3\right]_{q^2}\left[5\right]_{q^2}} 
{\left[10\right]_{q^2}\left[12\right]_{q^2}\left[15\right]_{q^2}} 
,\end{multline*} 
which, by Corollary~\ref{cor:A} and Lemma~\ref{lem:101215}, 
is a polynomial in $q$ with 
non-negative integer coefficients.

If $m\equiv
8
~(\text{mod }84),$ then we have \begin{multline*} \Cat^m(E_8;q)=
{\left[ 15 m+1 \right]_{q^{ 2 }}} {\left[ 
  \tfrac{15m+4}{4}  \right]_{q^{ 8 }}} {\left[ 
   \tfrac{5m+2}{42}  \right]_{q^{252}}} \frac {{ {\left[ 126
    \right]_{q^{ 2}}}{\left[ 3
    \right]_{q^{ 2}}} }} { {\left[ 6
    \right]_{q^{ 2}}}{\left[ 7 \right]_{q^{ 2}}} 
  {\left[ 9 \right]_{q^{ 2}}}}
 {\left[ 15 m+7 \right]_{q^{ 
  2 }}} {\left[ 5 m+3 \right]_{q^6}} \\
\times
 {\left[ \tfrac{3m+2}{2}  \right]_{q^{ 20 }}} 
 {\left[ \tfrac{5m+4}{4}  \right]_{q^{ 24}}}
{\left[ m+1 
  \right]_{q^{ 30 }}}
,\end{multline*} 
which, by Lemma~\ref{lem:679}, is a polynomial in $q$ with 
non-negative integer coefficients.

If $m\equiv
11
~(\text{mod }84),$ then we have 
\begin{multline*} \Cat^m(E_8;q)=
{\left[ \tfrac{15 m+1}2 \right]_{q^{ 4}}} 
{\left[ {15m+4} \right]_{q^{ 
   2}}} 
 {\left[ \tfrac{5m+2}3 \right]_{q^{ 18 }}}
 {\left[ \tfrac{15m+7}4 
   \right]_{q^{ 8}}} 
 {\left[ \tfrac{5m+3}2 \right]_{q^{ 12}}}\\
\times 
 {\left[ \tfrac{3m+2}7 \right]_{q^{ 70 }}} 
\frac {[35]_{q^2}} {[5]_{q^2}\left[7\right]_{q^2}}
{\left[ 5m+4
   \right]_{q^{ 6 }}} {\left[ \tfrac{m+1}4 \right]_{q^{ 120}}}
\frac {[60]_{q^2}\left[2\right]_{q^2}\left[3\right]_{q^2}\left[5\right]_{q^2}} 
{\left[10\right]_{q^2}\left[12\right]_{q^2}\left[15\right]_{q^2}} 
,\end{multline*} 
which, by Corollary~\ref{cor:A} and Lemma~\ref{lem:101215B}, 
is a polynomial in $q$ with 
non-negative integer coefficients.

If $m\equiv
16
~(\text{mod }84),$ then we have \begin{multline*} \Cat^m(E_8;q)=
{\left[ {15m+1}  \right]_{q^{ 2 }}}
 {\left[ \tfrac{15m+4}4 \right]_{q^{ 8 }}} {\left[ 
  \tfrac{5m+2}2  \right]_{q^{ 12 }}}
 {\left[ {15 m+7} 
  \right]_{q^2}} {\left[ {5m+3}  
   \right]_{q^{ 6}}}
   \\
\times {\left[ \tfrac{3m+2}2 \right]_{q^{20}}} 
 {\left[ \tfrac{5m+4}{84}  \right]_{q^{ 504 }}}
\frac {[252]_{q^2}\left[3\right]_{q^2}} 
{[7]_{q^2}\left[9\right]_{q^2}\left[12\right]_{q^2}}
 {\left[ {m+1} \right]_{q^{ 30}}} 
,\end{multline*} 
which, by Lemma~\ref{lem:7912}, 
is a polynomial in $q$ with 
non-negative integer coefficients.

If $m\equiv
18
~(\text{mod }84),$ then we have \begin{multline*} \Cat^m(E_8;q)=
{\left[ 15 m+1 \right]_{q^{ 2 }}} {\left[ \tfrac{15m+4}2 
  \right]_{q^4}} {\left[ \tfrac{5m+2}{4}  
   \right]_{q^{ 24}}} 
{\left[ 15 m+7 \right]_{q^2}} 
 {\left[ \tfrac{5m+3}{3}  \right]_{q^{ 18 }}} \\
 {\left[ \tfrac{3m+2}{28}  \right]_{q^{ 280}}} 
\frac {{ 
  {\left[ 140 \right]_{q^{ 2}}}  {\left[ 2 \right]_{q^{ 2}}} }}
 {  {\left[ 4 \right]_{q^{ 2}}} {\left[ 7 
   \right]_{q^{ 2}}} {\left[ 10 \right]_{q^{ 2}}}} {\left[ 
  \tfrac{5m+4}{2}  \right]_{q^{ 12 }}} {\left[ m+1 
  \right]_{q^{ 30 }}}
,\end{multline*} 
which, by Lemma~\ref{lem:4710}, 
is a polynomial in $q$ with 
non-negative integer coefficients.

If $m\equiv
21
~(\text{mod }84),$ then we have \begin{multline*} \Cat^m(E_8;q)=
{\left[ \tfrac{15m+1}{4}  \right]_{q^{ 8 }}}
 {\left[ 15m+4 \right]_{q^{ 2 }}} {\left[ 
  {5m+2}  \right]_{q^{ 6 }}}
 {\left[ \tfrac{15 m+7}{14} 
  \right]_{q^{28}}}
\frac {[14]_{q^2}} {[2]_{q^2}\left[7\right]_{q^2}}
 {\left[ \tfrac{5m+3}{12}  
   \right]_{q^{ 72}}}
\frac {[12]_{q^6}} {[3]_{q^6}\left[4\right]_{q^6}} 
   \\
\times {\left[ 3m+2 \right]_{q^{10}}} 
 {\left[ {5m+4}  \right]_{q^{ 6 }}}
 {\left[ \tfrac{m+1}{2} \right]_{q^{ 60}}} \frac {{ 
  {\left[ 30 \right]_{q^{ 2}}}
  {\left[ 2 \right]_{q^{ 2}}}
  {\left[ 3 \right]_{q^{ 2}}}
  {\left[ 5 \right]_{q^{ 2}}} }} { {\left[ 6 
   \right]_{q^{ 2}}}{\left[ 10
   \right]_{q^{ 2}}} {\left[ 15 \right]_{q^{ 2}}}}
,\end{multline*} 
which, by Corollary~\ref{cor:A} and Lemma~\ref{lem:61015B}, 
is a polynomial in $q$ with 
non-negative integer coefficients.

If $m\equiv
25
~(\text{mod }84),$ then we have \begin{multline*} \Cat^m(E_8;q)=
{\left[ \tfrac{15m+1}{4}  \right]_{q^{ 8 }}}
 {\left[ 15m+4 \right]_{q^{ 2 }}} {\left[ 
  {5m+2}  \right]_{q^{ 6 }}}
 {\left[ \tfrac{15 m+7}2 
  \right]_{q^4}} {\left[ \tfrac{5m+3}{4}  
   \right]_{q^{ 24}}} 
   \\
\times {\left[ \tfrac{3m+2}7 \right]_{q^{70}}} 
\frac {[35]_{q^2}} {[5]_{q^2}\left[7\right]_{q^2}}
 {\left[ \tfrac{5m+4}3  \right]_{q^{ 18 }}}
 {\left[ \tfrac{m+1}{2} \right]_{q^{ 60}}} \frac {{ 
  {\left[ 30 \right]_{q^{ 2}}}
  {\left[ 2 \right]_{q^{ 2}}}
  {\left[ 3 \right]_{q^{ 2}}}
  {\left[ 5 \right]_{q^{ 2}}} }} { {\left[ 6 
   \right]_{q^{ 2}}}{\left[ 10
   \right]_{q^{ 2}}} {\left[ 15 \right]_{q^{ 2}}}}
,\end{multline*} 
which, by Lemma~\ref{lem:61015C}, 
is a polynomial in $q$ with 
non-negative integer coefficients.

If $m\equiv
27
~(\text{mod }84),$ then we have 
\begin{multline*} \Cat^m(E_8;q)=
{\left[ \tfrac{15 m+1}{14} \right]_{q^{ 28}}}
\frac {[14]_{q^2}} {[2]_{q^2}\left[7\right]_{q^2}} 
{\left[ {15m+4} \right]_{q^{ 
   2}}} 
 {\left[ 5m+2 \right]_{q^{ 6 }}}
 {\left[ \tfrac{15m+7}4 
   \right]_{q^{ 8}}} 
 {\left[ \tfrac{5m+3}6 \right]_{q^{ 36}}}
\frac {[6]_{q^6}} {[2]_{q^6}\left[3\right]_{q^6}}
\\
\times 
 {\left[ 3m+2 \right]_{q^{ 10 }}} {\left[ 5m+4
   \right]_{q^{ 6 }}} {\left[ \tfrac{m+1}4 \right]_{q^{ 120}}}
\frac {[60]_{q^2}\left[2\right]_{q^2}\left[3\right]_{q^2}\left[5\right]_{q^2}} 
{\left[10\right]_{q^2}\left[12\right]_{q^2}\left[15\right]_{q^2}} 
,\end{multline*} 
which, by Corollary~\ref{cor:A} and Lemma~\ref{lem:101215C}, 
is a polynomial in $q$ with 
non-negative integer coefficients.

\medskip
All other cases are disposed of in a similar fashion.
\end{proof}

\section{Auxiliary results I}
\label{sec:aux1}

This section collects several auxiliary results which allow us to
reduce the problem of proving Theorem~\ref{thm:1}, or the
equivalent statement
\eqref{eq:1}, for the 26 exceptional groups listed in
Section~\ref{sec:prel} to a finite problem. While Lemmas~\ref{lem:2}
and \ref{lem:3} cover special choices of the parameters, 
Lemmas~\ref{lem:1} and
\ref{lem:6} afford an inductive procedure. More precisely, 
if we assume that we have already verified Theorem~\ref{thm:1} for all 
groups of smaller rank, then Lemmas~\ref{lem:1} and \ref{lem:6}, together 
with Lemmas~\ref{lem:2} and \ref{lem:7}, 
reduce the verification of Theorem~\ref{thm:1}
for the group that we are currently considering to a finite problem;
see Remark~\ref{rem:1}.
The final lemma of this section, Lemma~\ref{lem:8}, disposes of 
complex reflection groups with a special property satisfied by their degrees.

Let $p=am+b$, $0\le b<m$. We have
\begin{multline}
\phi^p\big((w_0;w_1,\dots,w_m)\big)\\
=(*;
c^{a+1}w_{m-b+1}c^{-a-1},c^{a+1}w_{m-b+2}c^{-a-1},
\dots,c^{a+1}w_{m}c^{-a-1},\\
c^{a}w_{1}c^{-a},\dots,
c^{a}w_{m-b}c^{-a}\big),
\label{eq:Aktion}
\end{multline}
where $*$ stands for the element of $W$ which is needed to 
complete the product of the components to $c$.

\begin{lemma} \label{lem:1}
It suffices to check \eqref{eq:1} for $p$ a divisor of $mh$.
More precisely, let $p$ be a divisor of $mh$, and let $k$ be another positive integer with 
$\gcd(k,mh/p)=1$, then we have
\begin{equation} \label{eq:2}
\Cat^m(W;q)\big\vert_{q=e^{2\pi i p/mh}}
= 
\Cat^m(W;q)\big\vert_{q=e^{2\pi i kp/mh}}
\end{equation}
and
\begin{equation} \label{eq:3}
\vert\Fix_{NC^m(W)}(\phi^{p})\vert =
\vert\Fix_{NC^m(W)}(\phi^{kp})\vert 	.
\end{equation}
\end{lemma}

\begin{proof}
For \eqref{eq:2}, this follows immediately from
\begin{equation} \label{eq:limit}
\lim_{q\to\zeta} \frac {[\alpha]_q} {[\beta]_q}=
\begin{cases}
\frac \alpha \beta&\text{if }\alpha\equiv\beta\equiv0\pmod d,\\
1&\text{otherwise},
\end{cases}	
\end{equation}
where $\zeta$ is a $d$-th root of unity and $\alpha,\beta$ are
non-negative integers such that 
$\alpha\equiv\beta\pmod d$.

In order to establish \eqref{eq:3}, suppose that 
$x\in\Fix_{NC^m(W)}(\phi^{p})$, that is, $x\in NC^m(W)$ and $\phi^p(x)=x$.
It obviously follows that $\phi^{kp}(x)=x$, so that
$x\in\Fix_{NC^m(W)}(\phi^{kp})$. To establish the converse, note that, 
if $\gcd(k,mh/p)=1$, then there exists $k'$ with $k'k\equiv 
1$~(mod~$\frac {mh} p$). It follows that, if 
$x\in\Fix_{NC^m(W)}(\phi^{kp})$, that is, if $x\in NC^m(W)$ and
$\phi^{kp}(x)=x$, then
$x=\phi^{k'kp}(x)=\phi^{p}(x)$, whence
$x\in\Fix_{NC^m(W)}(\phi^{p})$. 
\end{proof}

\begin{lemma} \label{lem:2}
Let $p$ be a divisor of $mh$.
If $p$ is divisible by $m$, then \eqref{eq:1} is true.
\end{lemma}


\begin{proof}
According to \eqref{eq:Aktion}, the action of $\phi^p$ on
$NC^m(W)$ is described by
\begin{equation*}
\phi^p\big((w_0;w_1,\dots,w_m)\big)
=(*;c^{p/m}w_{1}c^{-p/m},\dots,
c^{p/m}w_{m}c^{-p/m}\big).
\end{equation*}
Hence, if $(w_0;w_1,\dots,w_m)$ is fixed by $\phi^p$, then
each individual $w_i$ must be fixed under conjugation by $c^{p/m}$.

Using the notation $W'=\Cent_W(c^{p/m})$, 
the previous observation means that $w_i\in W'$, 
$i=1,2,\dots,m$. Springer \cite[Theorem~4.2]{SpriAA} 
(see also \cite[Theorem~11.24(iii)]{LeTaAA}) proved that
$W'$ is a well-generated complex reflection group whose degrees
coincide with those degrees of $W$ that are divisible by $mh/p$.
It was furthermore shown in \cite[Lemma~3.3]{BeReAA} that
\begin{equation} \label{eq:7}
NC(W)\cap W'=NC(W').
\end{equation}
Hence, the tuples $(w_0;w_1,\dots,w_m)$ fixed by $\phi^p$
are in fact identical with the elements of $NC^m(W')$, which
implies that
\begin{equation} \label{eq:6}
\vert\Fix_{NC^m(W)}(\phi^{p})\vert=\vert NC^m(W')\vert.
\end{equation}
Application of
Theorem~\ref{thm:2} with $W$ replaced by $W'$ and of the
``limit rule" \eqref{eq:limit} then yields that
\begin{equation} \label{eq:5}
\vert NC^m(W')\vert=
\underset{\frac {mh} p\mid d_i}{\prod_{1\le i\le n}} \frac {mh+d_i} {d_i}=\Cat^m(W;q)\big\vert_{q=e^{2\pi i p/mh}}.
\end{equation}
Combining \eqref{eq:6} and \eqref{eq:5}, we obtain \eqref{eq:1}.
This finishes the proof of the lemma.
\end{proof}

\begin{lemma} \label{lem:3}
Equation \eqref{eq:1} holds for all divisors $p$ of $m$.
\end{lemma}

\begin{proof}
Using \eqref{eq:limit} and the fact that the degrees of 
irreducible well-generated complex reflection groups satisfy
$d_i<h$ for all $i<n$, we see that
$$
\Cat^m(W;q)\big\vert_{q=e^{2\pi i p/mh}}=\begin{cases}
m+1&\text{if }m=p,\\
1&\text{if }m\ne p.
\end{cases}
$$
On the other hand, if $(w_0;w_1,\dots,w_m)$ is fixed by $\phi^p$,
then, because of the action \eqref{eq:Aktion}, we must have
$w_1=w_{p+1}=\dots=w_{m-p+1}$ and $w_1=cw_{m-p+1}c^{-1}$. In particular,
$w_1\in\Cent_W(c)$. By the theorem of Springer cited in the
proof of Lemma~\ref{lem:2}, the subgroup $\Cent_W(c)$
is itself a complex reflection group whose degrees are those degrees
of $W$ that are divisible by $h$. The only such degree is $h$
itself, hence $\Cent_W(c)$ is the cyclic group generated by
$c$. Moreover, by \eqref{eq:7}, we obtain that $w_1=\ep$,
the identity element of $W$, or
$w_1=c$. Therefore, for $m=p$ the set $\Fix_{NC^m(W)}(\phi^p)$
consists of the $m+1$ elements $(w_0;w_1,\dots,w_m)$ obtained by choosing 
$w_i=c$ for a particular $i$ between $0$ and $m$, all other $w_j$'s
being equal to $\ep$, while, for $m\ne p$, we have
$$\Fix_{NC^m(W)}(\phi^p)=\big\{(c;\ep,\dots,\ep)\big\},$$
whence the result.
\end{proof}

\begin{lemma} \label{lem:4}
Let $W$ be an irreducible well-generated complex reflection group 
all of whose degrees are divisible by $d$. 
Then each element of $W$ is fixed under conjugation by $c^{h/d}$.
\end{lemma}

\begin{proof}
By the theorem of Springer cited in the
proof of Lemma~\ref{lem:2}, the subgroup $W'=\Cent_W(c^{h/d})$
is itself a complex reflection group whose degrees are those degrees
of $W$ that are divisible by $d$. Thus, by our assumption, the
degrees of $W'$ coincide with the degrees of $W$, and hence $W'$ must
be equal to $W$.
Phrased differently, each element of $W$ is fixed under conjugation
by $c^{h/d}$, as claimed.
\end{proof}

\begin{lemma} \label{lem:6}
Let $W$ be an irreducible 
well-generated complex reflection group of rank $n$, 
and let $p=m_1h_1$ be a divisor of $mh$, where $m=m_1m_2$ and
$h=h_1h_2$. Without loss of generality, we assume that 
$\gcd(h_1,m_2)=1$. Suppose that
Theorem~{\em\ref{thm:1}} has already been verified for all 
irreducible well-generated
complex reflection groups with rank $<n$. 
If $h_2$ does not divide all degrees $d_i$,
then Equation~\eqref{eq:1} is satisfied.
\end{lemma}

\begin{proof}
Let us write $h_1=am_2+b$, with $0\le b<m_2$. The condition 
$\gcd(h_1,m_2)=1$ translates into $\gcd(b,m_2)=1$.
From \eqref{eq:Aktion}, we infer that
\begin{multline} \label{eq:m2Aktion}
\phi^p\big((w_0;w_1,\dots,w_m)\big)\\=
(*;
c^{a+1}w_{m-m_1b+1}c^{-a-1},c^{a+1}w_{m-m_1b+2}c^{-a-1},
\dots,c^{a+1}w_{m}c^{-a-1},\\
c^aw_{1}c^{-a},\dots,
c^aw_{m-m_1b}c^{-a}\big).
\end{multline}
Supposing that 
$(w_0;w_1,\dots,w_m)$ is fixed by $\phi^p$, we obtain
the system of equations
\begin{align*} 
w_i&=c^{a+1}w_{i+m-m_1b}c^{-a-1}, \quad i=1,2,\dots,m_1b,\\
w_i&=c^aw_{i-m_1b}c^{-a}, \quad i=m_1b+1,m_1b+2,\dots,m,
\end{align*}
which, after iteration, implies in particular that
$$
w_i=c^{b(a+1)+(m_2-b)a}w_ic^{-b(a+1)-(m_2-b)a}=c^{h_1}w_ic^{-h_1},
\quad i=1,2,\dots,m.
$$
It is at this point where we need $\gcd(b,m_2)=1$.
The last equation shows that each $w_i$, $i=1,2,\dots,m$, and thus
also $w_0$, lies in $\Cent_{W}(c^{h_1})$. 
By the theorem of Springer cited in the
proof of Lemma~\ref{lem:2}, this centraliser subgroup 
is itself a complex reflection group, $W'$ say, 
whose degrees are those degrees
of $W$ that are divisible by $h/h_1=h_2$. Since, by assumption, $h_2$
does not divide {\em all\/} degrees, $W'$ has 
rank strictly less than $n$. Again by assumption, we know that
Theorem~\ref{thm:1} is true for $W'$, so that in particular,
$$
\vert\Fix_{NC^m(W')}(\phi^{p})\vert = 
\Cat^m(W';q)\big\vert_{q=e^{2\pi i p/mh}}.
$$
The arguments above together with \eqref{eq:7} show that 
$\Fix_{NC^m(W)}(\phi^{p})=\Fix_{NC^m(W')}(\phi^{p})$.
On the other hand, using \eqref{eq:limit} 
it is straightforward to see that
$$\Cat^m(W;q)\big\vert_{q=e^{2\pi i p/mh}}=
\Cat^m(W';q)\big\vert_{q=e^{2\pi i p/mh}}.$$
This proves \eqref{eq:1} for our particular $p$, as required. 
\end{proof}

\begin{lemma} \label{lem:7}
Let $W$ be an irreducible 
well-generated complex reflection group of rank $n$, 
and let $p=m_1h_1$ be a divisor of $mh$, where $m=m_1m_2$ and
$h=h_1h_2$. We assume that 
$\gcd(h_1,m_2)=1$. If $m_2>n$ then
$$
\Fix_{NC^m(W)}(\phi^{p})=\big\{(c;\ep,\dots,\ep)\big\}.
$$
\end{lemma}

\begin{proof}
Let us suppose that $(w_0;w_1,\dots,w_m)\in 
\Fix_{NC^m(W)}(\phi^{p})$ and that there exists a $j\ge1$ such
that $w_j\ne\ep$. By \eqref{eq:m2Aktion}, it then follows
for such a $j$ that
also $w_k\ne\ep$ for all $k\equiv j-lm_1b$~(mod~$m$), where, as before,
$b$ is defined as the unique integer with $h_1=am_2+b$ and
$0\le b<m_2$. Since, by assumption, $\gcd(b,m_2)=1$, there are
exactly $m_2$ such $k$'s which are distinct mod~$m$. 
However, this implies that the sum of the absolute lengths
of the $w_i$'s, $0\le i\le m$, is at least $m_2>n$, a
contradiction to Remark~\ref{rem:0}.(2).
\end{proof}

\begin{remark} \label{rem:1}
(1)
If we put ourselves in the situation of the assumptions of
Lemma~\ref{lem:6}, then we may conclude that equation~\eqref{eq:1} 
only needs to be checked for pairs $(m_2,h_2)$ subject to the 
following restrictions: 
\begin{equation} 
m_2\ge2,\quad 
\gcd(h_1,m_2)=1,\quad
\text{and $h_2$ divides all degrees of $W$}. 
\label{eq:restr}
\end{equation}
Indeed, Lemmas~\ref{lem:2} and \ref{lem:6} 
together
imply that equation~\eqref{eq:1} is always satisfied 
in all other cases.

\smallskip
(2) Still putting ourselves in the situation of Lemma~\ref{lem:6},
if $m_2>n$ and $m_2h_2$ does not divide any of the degrees of $W$,
then equation~\eqref{eq:1} is satisfied. Indeed, Lemma~\ref{lem:7}
says that in this case the left-hand side of \eqref{eq:1} equals
$1$, while a straightforward computation using \eqref{eq:limit}
shows that in this case the right-hand side of \eqref{eq:1}
equals $1$ as well. 

\smallskip
(3) It should be observed that this leaves a finite
number of choices for $m_2$ to consider, whence a finite number of
choices for $(m_1,m_2,h_1,h_2)$. Altogether, there remains a finite
number of choices for $p=h_1m_1$ to be checked.
\end{remark}

\begin{lemma} \label{lem:8}
Let $W$ be an irreducible 
well-generated complex reflection group of rank $n$ 
with the property that $d_i\mid h$ for $i=1,2,\dots,n$.
Then Theorem~{\em\ref{thm:1}} is true for this group $W$.
\end{lemma}

\begin{proof}
By Lemma~\ref{lem:1}, we may restrict ourselves to divisors
$p$ of $mh$. 

Suppose that $e^{2\pi ip/mh}$ is a $d_i$-th root of unity
for some $i$. In other words, $mh/p$ divides $d_i$.
Since $d_i$ is a divisor of $h$ by assumption, 
the integer $mh/p$ also divides $h$. But this is equivalent to
saying that $m$ divides $p$, and equation \eqref{eq:1} holds by 
Lemma~\ref{lem:2}.

Now assume that $mh/p$ does not divide any of the $d_i$'s.
Then, by \eqref{eq:limit}, the right-hand side of \eqref{eq:1}
equals $1$.  On the other hand, $(c;\ep,\dots,\ep)$ is always an
element of $\Fix_{NC^m(W)}(\phi^{p})$. To see that there are no
others, we make appeal to the classification of all irreducible 
well-generated complex reflection groups, which we recalled in
Section~\ref{sec:prel}. Inspection reveals that all groups
satisfying the hypotheses of the lemma have rank $n\le2$.
Except for the groups contained in the infinite series $G(d,1,n)$ 
and $G(e,e,n)$ for which Theorem~\ref{thm:1} has been established in
\cite{KratCG}, these are the groups $G_5,G_6,G_9,G_{10},
G_{14},G_{17},G_{18},G_{21}$. We now discuss these groups case
by case, keeping the notation of Lemma~\ref{lem:6}.
In order to simplify the argument, we note that Lemma~\ref{lem:7}
implies that equation~\eqref{eq:1} holds if $m_2>2$, so that
in the following arguments we always may assume that $m_2=2$.


\smallskip
{\sc Case $G_5$}. The degrees are $6,12$, and
therefore Remark~\ref{rem:1}.(1) implies that equation~\eqref{eq:1}
is always satisfied.

\smallskip
{\sc Case $G_6$}. The degrees are $4,12$, and
therefore, according to Remark~\ref{rem:1}.(1), we need only consider
the case where $h_2=4$ and $m_2=2$, that is, $p=3m/2$. Then \eqref{eq:m2Aktion} becomes
\begin{equation} \label{eq:3m2Aktion}
\phi^p\big((w_0;w_1,\dots,w_m)\big)
=(*;
c^{2}w_{\frac m2+1}c^{-2},
c^{2}w_{\frac m2+2}c^{-2},
\dots,c^{2}w_{m}c^{-2},
cw_{1}c^{-1},\dots,
cw_{\frac m2}c^{-1}\big).
\end{equation}
If $(w_0;w_1,\dots,w_m)$ is fixed by $\phi^p$ and not
equal to $(c;\ep,\dots,\ep)$, there must exist an $i$
with $1\le i\le \frac m2$ such that 
$\ell_T(w_i)=\ell_T(w_{\frac {m} {2}+i})=1$,
$w_{\frac {m} {2}+i}=cw_ic^{-1}$,
$w_iw_{\frac {m} {2}+i}=w_icw_{i}c^{-1}=c$, and all $w_j$, 
with $j\ne i,\frac m2+i$,
equal $\ep$. However, 
with the help of the {\sl GAP} package {\tt CHEVIE}
\cite{chevAA,MichAA}, one verifies that there is no $w_i$
in $G_6$ such that
$$
\ell_T(w_i)=1\quad \text{and}\quad 
w_icw_{i}c^{-1}=c
$$
are simultaneously satisfied. Hence,
the left-hand side of \eqref{eq:1} is equal to $1$, as required.

\smallskip
{\sc Case $G_9$}. The degrees are $8,24$, and
therefore, according to Remark~\ref{rem:1}.(1), we need only consider
the case where $h_2=8$ and $m_2=2$, that is, $p=3m/2$. 
This is the same $p$ as for $G_6$. Again, 
{\tt CHEVIE}
finds no solution. Hence,
the left-hand side of \eqref{eq:1} is equal to $1$, as required.

\smallskip
{\sc Case $G_{10}$}. The degrees are $12,24$, and
therefore Remark~\ref{rem:1}.(1) implies that equation~\eqref{eq:1}
is always satisfied.

\smallskip
{\sc Case $G_{14}$}. The degrees are $6,24$, and
therefore Remark~\ref{rem:1}.(1) implies that equation~\eqref{eq:1}
is always satisfied.

\smallskip
{\sc Case $G_{17}$}. The degrees are $20,60$, and
therefore, according to Remark~\ref{rem:1}.(1), we need only consider
the cases where $h_2=20$ or $h_2=4$. 
In the first case, $p=3m/2$, which is
the same $p$ as for $G_6$. Again, 
{\tt CHEVIE} finds no solution. 
In the second case, $p=15m/2$. Then \eqref{eq:m2Aktion} becomes
\begin{multline} \label{eq:15m2Aktion}
\phi^p\big((w_0;w_1,\dots,w_m)\big)\\
=(*;
c^{8}w_{\frac m2+1}c^{-8},
c^{8}w_{\frac m2+2}c^{-8},
\dots,c^{8}w_{m}c^{-8},
c^7w_{1}c^{-7},\dots,
c^7w_{\frac m2}c^{-7}\big).
\end{multline}
By Lemma~\ref{lem:4}, every element of $NC(W)$ is fixed under 
conjugation by $c^3$, and, thus, on elements fixed by $\phi^p$,
the above action of $\phi^p$
reduces to the one in \eqref{eq:3m2Aktion}. This action was already
discussed in the first case.
Hence, in both cases,
the left-hand side of \eqref{eq:1} is equal to $1$, as required.

\smallskip
{\sc Case $G_{18}$}. The degrees are $30,60$, and
therefore Remark~\ref{rem:1}.(1) implies that equation~\eqref{eq:1}
is always satisfied.

\smallskip
{\sc Case $G_{21}$}. The degrees are $12,60$, and
therefore, according to Remark~\ref{rem:1}.(1), we need only consider
the cases where $h_2=12$ or $h_2=4$. 
In the first case, $p=5m/2$, 
so that \eqref{eq:m2Aktion} becomes
\begin{multline} \label{eq:5m2Aktion}
\phi^p\big((w_0;w_1,\dots,w_m)\big)\\
=(*;
c^{3}w_{\frac m2+1}c^{-3},
c^{3}w_{\frac m2+2}c^{-3},
\dots,c^{3}w_{m}c^{-3},
c^2w_{1}c^{-2},\dots,
c^2w_{\frac m2}c^{-2}\big).
\end{multline}
If $(w_0;w_1,\dots,w_m)$ is fixed by $\phi^p$ and not
equal to $(c;\ep,\dots,\ep)$, there must exist an $i$
with $1\le i\le \frac m2$ such that $\ell_T(w_i)=1$ and
$w_ic^2w_{i}c^{-2}=c$. However, 
with the help of the {\sl GAP} package {\tt CHEVIE}
\cite{chevAA,MichAA}, one verifies that there is no such
solution to this equation. 
In the second case, $p=15m/2$. Then \eqref{eq:m2Aktion} becomes
the action in \eqref{eq:15m2Aktion}.
By Lemma~\ref{lem:4}, every element of $NC(W)$ is fixed under 
conjugation by $c^5$, and, thus, on elements fixed by $\phi^p$,
the action of $\phi^p$ in
\eqref{eq:15m2Aktion}
reduces to the one in the first case. 
Hence, in both cases,
the left-hand side of \eqref{eq:1} is equal to $1$, as required.

\smallskip
This completes the proof of the lemma.
\end{proof}

\section{Exemplification of case-by-case verification of Theorem~\ref{thm:1}}
\label{sec:Beweis1}

It remains to verify Theorem~\ref{thm:1} for the groups
$G_4,G_8,G_{16},G_{20},
G_{23}=H_3,G_{24},G_{25},\break G_{26},G_{27},
G_{28}=F_4,G_{29},G_{30}=H_4,G_{32},
G_{33},G_{34},G_{35}=E_6,G_{36}=E_7,G_{37}=E_8$.
All details can be found in \cite[Sec.~6]{KrMuAE}.
We content ourselves with illustrating the 
type of computation that is needed here by going through the case
of the group $G_{24}$, and by discussing some of the arguments needed
for the group $G_{37}=E_8$. 

In the sequel we write $\zeta_d$ for a primitive $d$-th root of 
unity.

\subsection*{\sc Case $G_{24}$}
The degrees are $4,6,14$, and hence we have
$$
\Cat^m(G_{24};q)=\frac 
{[14m+14]_q\, [14m+6]_q\, [14m+4]_q} 
{[14]_q\, [6]_q\, [4]_q} .
$$
Let $\zeta$ be a $14m$-th root of unity. 
In what follows, 
we abbreviate the assertion that ``$\zeta$ is a primitive $d$-th root of
unity" as ``$\zeta=\zeta_d$."
The following cases on the right-hand side of \eqref{eq:1}
occur:
{\refstepcounter{equation}\label{eq:G24}}
\alphaeqn
\begin{align} 
\label{eq:G24.2}
\lim_{q\to\zeta}\Cat^m(G_{24};q)&=m+1,
\quad\text{if }\zeta=\zeta_{14},\zeta_7,\\
\label{eq:G24.3}
\lim_{q\to\zeta}\Cat^m(G_{24};q)&=\tfrac {7m+3}3,
\quad\text{if }\zeta=\zeta_{6},\zeta_3,\ 3\mid m,\\
\label{eq:G24.4}
\lim_{q\to\zeta}\Cat^m(G_{24};q)&=\tfrac {7m+2}2,
\quad\text{if }\zeta=\zeta_4,\ 2\mid m,\\
\label{eq:G24.5}
\lim_{q\to\zeta}\Cat^m(G_{24};q)&=\Cat^m(G_{24}),
\quad\text{if }\zeta=-1\text{ or }\zeta=1,\\
\label{eq:G24.1}
\lim_{q\to\zeta}\Cat^m(G_{24};q)&=1,
\quad\text{otherwise.}
\end{align}
\reseteqn

We must now prove that the left-hand side of \eqref{eq:1} in
each case agrees with the values exhibited in 
\eqref{eq:G24}. The only cases not covered by
Lemma~\ref{lem:2} are the ones in \eqref{eq:G24.3},
\eqref{eq:G24.4},
and \eqref{eq:G24.1}. (In both \eqref{eq:G24.2} and \eqref{eq:G24.5}
we have $d\mid h$.)

We first consider \eqref{eq:G24.3}. 
By Lemma~\ref{lem:1}, we are free to choose $p=7m/3$ if 
$\zeta=\zeta_6$, respectively $p=14m/3$ if 
$\zeta=\zeta_3$. In both cases, $m$ must be divisible by $3$.

We start with the case that $p=7m/3$.
From \eqref{eq:Aktion}, we infer
\begin{multline*}
\phi^p\big((w_0;w_1,\dots,w_m)\big)\\=
(*;
c^{3}w_{\frac {2m}3+1}c^{-3},c^{3}w_{\frac {2m}3+2}c^{-3},
\dots,c^{3}w_{m}c^{-3},
c^2w_{1}c^{-2},\dots,
c^2w_{\frac {2m}3}c^{-2}\big).
\end{multline*}
Supposing that 
$(w_0;w_1,\dots,w_m)$ is fixed by $\phi^p$, we obtain
the system of equations
{\refstepcounter{equation}\label{eq:G24A}}
\alphaeqn
\begin{align} \label{eq:G24Aa}
w_i&=c^3w_{\frac {2m}3+i}c^{-3}, \quad i=1,2,\dots,\tfrac {m}3,\\
w_i&=c^2w_{i-\frac {m}3}c^{-2}, \quad i=\tfrac {m}3+1,\tfrac {m}3+2,\dots,m.
\label{eq:G24Ab}
\end{align}
\reseteqn
There are two distinct possibilities for choosing
the $w_i$'s, $1\le i\le m$: 
either all the $w_i$'s are equal to $\ep$, or
there is an $i$ with $1\le i\le \frac m3$ such that
$$\ell_T(w_i)=\ell_T(w_{i+\frac m3})=\ell_T(w_{i+\frac {2m}3})=1.$$
Writing $t_1,t_2,t_3$ for $w_i,w_{i+\frac m3},w_{i+\frac {2m}3}$, 
respectively, the equations \eqref{eq:G24A} 
reduce to
{\refstepcounter{equation}\label{eq:G24B}}
\alphaeqn
\begin{align} \label{eq:G24Ba}
t_1&=c^3t_3c^{-3},\\
\label{eq:G24Bb}
t_2&=c^2t_1c^{-2},\\
t_3&=c^2t_2c^{-2}.
\label{eq:G24Bc}
\end{align}
\reseteqn
One of these equations is in fact superfluous: if we substitute
\eqref{eq:G24Bb} and \eqref{eq:G24Bc} in \eqref{eq:G24Ba}, then
we obtain $t_1=c^7t_1c^{-7}$ which is automatically satisfied due
to Lemma~\ref{lem:4} with $d=2$.

Since $(w_0;w_1,\dots,w_m)\in NC^m(G_{24})$, we must have $t_1t_2t_3=c$.
Combining this with \eqref{eq:G24B}, we infer that
\begin{equation} \label{eq:G24D}
t_1(c^{2}t_1c^{-2})(c^4t_1c^{-4})=c.
\end{equation}
With the help of {\tt CHEVIE}, 
one obtains 7 solutions for $t_1$ in this equation,
each of them giving rise to $m/3$ elements of
$\Fix_{NC^m(G_{24})}(\phi^{p})$ since $i$ (in $w_i$) ranges from $1$ to $m/3$.

In total, we obtain  
$1+7\frac m3=\frac {7m+3}3$ elements in
$\Fix_{NC^m(G_{24})}(\phi^p)$, which agrees with the limit in
\eqref{eq:G24.3}.

The case where $p=14m/3$ can be treated in a similar fashion.
In the end, it turns out that we have to solve the same enumeration 
problem as for
$p=7m/3$, and, consequently, the number of elements of
$\Fix_{NC^m(G_{24})}(\phi^{p})$ is the same, namely
$\frac {7m+3}3$, as required.

\smallskip
Our next case is \eqref{eq:G24.4}. 
Proceeding in a similar manner as before, we see that there is again the
trivial possibility $(c;\ep,\dots,\ep)$, and otherwise we have to
find $t_1$ with $\ell_T(t_1)=1$ satisfying the inequality
\begin{equation} \label{eq:G24''D}
t_1(c^{3}t_1c^{-3})\le_T c.
\end{equation}
With the help of {\tt CHEVIE}, 
one obtains 7 solutions for $t_1$ in this relation,
each of them giving rise to $m/2$ elements of
$\Fix_{NC^m(G_{24})}(\phi^{p})$ since $i$ (in $w_i$) ranges from $1$ to $m/2$.

In total, we obtain  
$1+7\frac m2=\frac {7m+2}2$ elements in
$\Fix_{NC^m(G_{24})}(\phi^p)$, which agrees with the limit in
\eqref{eq:G24.4}.

\smallskip
Finally, we turn to \eqref{eq:G24.1}. By Remark~\ref{rem:1},
the only choices for $h_2$ and $m_2$ to be considered
are $h_2=1$ and $m_2=3$, $h_2=m_2=2$, and $h_2=2$
and $m_2=3$. These correspond to the choices $p=14m/3$,
$p=7m/2$, respectively $p=7m/3$, all of which have already been 
discussed as they do not belong to \eqref{eq:G24.1}. Hence, 
\eqref{eq:1} must necessarily hold, as required.

\subsection*{\sc Case $G_{37}=E_8$}
The degrees are $2,8,12,14,18,20,24,30$, and hence we have
\begin{multline*}
\Cat^m(E_8;q)=\frac 
{[30m+30]_q\, [30m+24]_q\, [30m+20]_q\, [30m+18]_q} 
{[30]_q\, [24]_q\, [20]_q\, [18]_q}\\
\times
\frac 
{[30m+14]_q\, [30m+12]_q\, [30m+8]_q\, [30m+2]_q} 
{[14]_q\, [12]_q\, [8]_q\, [2]_q} .
\end{multline*}
Let $\zeta$ be a $30m$-th root of unity. 
The cases occurring on the right-hand side of \eqref{eq:1}
not covered by Lemma~\ref{lem:2} are:
{\refstepcounter{equation}\label{eq:E8}}
\alphaeqn
{\allowdisplaybreaks
\begin{align} 
\label{eq:E8.3}
\lim_{q\to\zeta}\Cat^m(E_8;q)&=\tfrac {5m+4}4,
\quad\text{if }\zeta=\zeta_{24},\ 4\mid m,\\
\label{eq:E8.4}
\lim_{q\to\zeta}\Cat^m(E_8;q)&=\tfrac {3m+2}2,
\quad\text{if }\zeta=\zeta_{20},\ 2\mid m,\\
\label{eq:E8.5}
\lim_{q\to\zeta}\Cat^m(E_8;q)&=\tfrac {5m+3}3,
\quad\text{if }\zeta=\zeta_{18},\zeta_{9},\ 3\mid m,\\
\label{eq:E8.6}
\lim_{q\to\zeta}\Cat^m(E_8;q)&=\tfrac {15m+7}7,
\quad\text{if }\zeta=\zeta_{14},\zeta_{7},\ 7\mid m,\\
\label{eq:E8.7}
\lim_{q\to\zeta}\Cat^m(E_8;q)&=\tfrac {(5m+4)(5m+2)}8,
\quad\text{if }\zeta= \zeta_{12},\ 2\mid m,\\
\label{eq:E8.9}
\lim_{q\to\zeta}\Cat^m(E_8;q)&=\frac{(5m+4)(15m+4)}{16},
\quad\text{if }\zeta= \zeta_{8},\ 4\mid m,\\
\label{eq:E8.11}
\lim_{q\to\zeta}\Cat^m(E_8;q)&=\frac {(5m+4)(3m+2)(5m+2)(15m+4)}{64},
\quad\text{if }\zeta= \zeta_{4},\ 2\mid m,\\
\lim_{q\to\zeta}\Cat^m(E_8;q)&=\Cat^m(E_8),
\quad\text{if }\zeta=-1\text{ or }\zeta=1,\\
\label{eq:E8.1}
\lim_{q\to\zeta}\Cat^m(E_8;q)&=1,
\quad\text{otherwise.}
\end{align}}%
\reseteqn

We now have to prove that the left-hand side of \eqref{eq:1} in
each case agrees with the values exhibited in 
\eqref{eq:E8}. Since the corresponding computations in the various 
cases are very similar, we concentrate here only on the cases
\eqref{eq:E8.9} and \eqref{eq:E8.11}, these two being representative
of the types of arguments arising. 
As before, we refer the reader to \cite[Sec.~6]{KrMuAE} for
full details. 

\smallskip
Let us consider the case in \eqref{eq:E8.9} first.
By Lemma~\ref{lem:1}, we are free to choose $p=15m/4$. In particular,
$m$ must be divisible by $4$.
From \eqref{eq:Aktion}, we infer
\begin{multline*}
\phi^p\big((w_0;w_1,\dots,w_m)\big)\\=
(*;
c^{4}w_{\frac {m}4+1}c^{-4},c^{4}w_{\frac {m}4+2}c^{-4},
\dots,c^{4}w_{m}c^{-4},
c^3w_{1}c^{-3},\dots,
c^3w_{\frac {m}4}c^{-3}\big).
\end{multline*}
Supposing that 
$(w_0;w_1,\dots,w_m)$ is fixed by $\phi^p$, we obtain
the system of equations
{\refstepcounter{equation}\label{eq:E8''A}}
\alphaeqn
\begin{align} \label{eq:E8''Aa}
w_i&=c^4w_{\frac {m}4+i}c^{-4}, \quad i=1,2,\dots,\tfrac {3m}4,\\
w_i&=c^3w_{i-\frac {3m}4}c^{-3}, \quad i=\tfrac {3m}4+1,\tfrac {3m}4+2,\dots,m.
\label{eq:E8''Ab}
\end{align}
\reseteqn
There are several distinct possibilities for choosing
the $w_i$'s, $1\le i\le m$, which we summarise as follows: 
{\refstepcounter{equation}\label{eq:E8''C}}
\alphaeqn
\begin{enumerate}
\item[(i)]
all the $w_i$'s are equal to $\ep$ (and $w_0=c$), 
\item[(ii)]
there is an $i$ with $1\le i\le \frac m4$ such that
\begin{equation} \label{eq:E8''Cii}
1\le\ell_T(w_i)=\ell_T(w_{i+\frac m4})=\ell_T(w_{i+\frac {2m}4})=
\ell_T(w_{i+\frac {3m}4})\le2,
\end{equation}
and the other $w_j$'s, $1\le j\le m$, are equal to $\ep$,
\item[(iii)]
there are $i_1$ and $i_2$ with $1\le i_1<i_2\le \frac m4$ such that
\begin{multline} \label{eq:E8''Ciii}
\ell_T(w_{i_1})=\ell_T(w_{i_2})=
\ell_T(w_{i_1+\frac m4})=\ell_T(w_{i_2+\frac m4})\\=
\ell_T(w_{i_1+\frac {2m}4})=\ell_T(w_{i_2+\frac {2m}4})=
\ell_T(w_{i_1+\frac {3m}4})=\ell_T(w_{i_2+\frac {3m}4})=1,
\end{multline}
and all other $w_j$ are equal to $\ep$.
\end{enumerate}
\reseteqn

Moreover, since $(w_0;w_1,\dots,w_m)\in NC^m(E_8)$, we must have 
$$w_iw_{i+\frac {m}4}w_{i+\frac {2m}4}w_{i+\frac {3m}4}\le_T c,$$
or
$$
w_{i_1}w_{i_2}w_{i_1+\frac {m}4}w_{i_2+\frac {m}4}
w_{i_1+\frac {2m}4}w_{i_2+\frac {2m}4}
w_{i_1+\frac {3m}4}w_{i_2+\frac {3m}4}=c.
$$
Together with equations~\eqref{eq:E8''A}--\eqref{eq:E8''C}, 
this implies that
\begin{equation} \label{eq:E8''D}
w_i=c^{15}w_ic^{-15}\quad\text{and}\quad 
w_i(c^{11}w_ic^{-11})(c^7w_ic^{-7})
(c^3w_ic^{-3})\le_T c,
\end{equation}
or that
\begin{multline} \label{eq:E8''E}
w_{i_1}=c^{15}w_{i_1}c^{-15},\quad
w_{i_1}=c^{15}w_{i_2}c^{-15},\\
\quad\text{and}\quad 
w_{i_1}w_{i_2}
(c^{11}w_{i_1}c^{-11})(c^{11}w_{i_2}c^{-11})
(c^7w_{i_1}c^{-7})(c^7w_{i_2}c^{-7})
(c^3w_{i_1}c^{-3})(c^3w_{i_2}c^{-3})=c.
\end{multline}
Here, the first equation in \eqref{eq:E8''D} and the first
two equations in \eqref{eq:E8''E} are 
automatically satisfied due to Lemma~\ref{lem:4} with $d=2$.

With the help of Stembridge's {\sl Maple} package {\tt coxeter}
\cite{StemAZ}, one obtains 30 solutions for $w_i$ in 
\eqref{eq:E8''D} with $\ell_T(w_i)=1$,
45 solutions for $w_i$ with $\ell_T(w_i)=2$ and $w_i$ of type $A_1^2$
(as a parabolic Coxeter element; see the end of Section~\ref{sec:prel}),
and 20 solutions for $w_i$ with $\ell_T(w_i)=2$ and $w_i$ of type $A_2$.
Each of them gives rise to $m/4$ elements of
$\Fix_{NC^m(E_8)}(\phi^{p})$ since $i$ ranges from $1$ to $m/4$.

The number of solutions in Case~(iii) can be computed
from our knowledge of the solutions in Case~(ii) according to type,
using some elementary counting arguments. 
Namely, the number of
solutions of \eqref{eq:E8''E} is equal to
$$
45\cdot 2+20\cdot 3=150,
$$
since an element of type $A_1^2$ can be decomposed in two ways 
into a product of two elements of absolute length $1$, while for
an element of type $A_2$ this can be done in $3$
ways. 

In total, we obtain 
$1+(30+45+20)\frac m4+150\binom {m/4}2=\frac {(5m+4)(15m+4)}{16}$ elements in
$\Fix_{NC^m(E_8)}(\phi^p)$, which agrees with the limit in
\eqref{eq:E8.9}.

\smallskip
Next, we discuss the case in \eqref{eq:E8.11}.
By Lemma~\ref{lem:1}, we are free to choose $p=15m/2$. In particular,
$m$ must be divisible by $2$.
From \eqref{eq:Aktion}, we infer
\begin{multline*}
\phi^p\big((w_0;w_1,\dots,w_m)\big)\\=
(*;
c^{8}w_{\frac {m}2+1}c^{-8},c^{8}w_{\frac {m}2+2}c^{-8},
\dots,c^{8}w_{m}c^{-8},
c^7w_{1}c^{-7},\dots,
c^7w_{\frac {m}2}c^{-7}\big).
\end{multline*}
Supposing that 
$(w_0;w_1,\dots,w_m)$ is fixed by $\phi^p$, we obtain
the system of equations
{\refstepcounter{equation}\label{eq:E8''''A}}
\alphaeqn
\begin{align} \label{eq:E8''''Aa}
w_i&=c^8w_{\frac {m}2+i}c^{-8}, \quad i=1,2,\dots,\tfrac {m}2,\\
w_i&=c^7w_{i-\frac {m}2}c^{-7}, \quad i=\tfrac {m}2+1,\tfrac {m}2+2,\dots,m.
\label{eq:E8''''Ab}
\end{align}
\reseteqn
There are several distinct possibilities for choosing
the $w_i$'s, $1\le i\le m$: 
{\refstepcounter{equation}\label{eq:E8''''C}}
\alphaeqn
\begin{enumerate}
\item[(i)]
all the $w_i$'s are equal to $\ep$ (and $w_0=c$), 
\item[(ii)]
there is an $i$ with $1\le i\le \frac m2$ such that
\begin{equation} \label{eq:E8''''Cii}
1\le\ell_T(w_i)=\ell_T(w_{i+\frac m2})\le 4,
\end{equation}
and the other $w_j$'s, $1\le j\le m$, are equal to $\ep$,
\item[(iii)]
there are $i_1$ and $i_2$ with $1\le i_1<i_2\le \frac m2$ such that
\begin{equation} \label{eq:E8''''Ciii}
\ell_1:=\ell_T(w_{i_1})=\ell_T(w_{i_1+\frac m2})\ge1,\quad
\ell_2:=\ell_T(w_{i_2})=\ell_T(w_{i_2+\frac m2})\ge1,\quad\text{and}\quad
\ell_1+\ell_2\le4,
\end{equation}
and the other $w_j$'s, $1\le j\le m$, are equal to $\ep$,
\item[(iv)]
there are $i_1,i_2,i_3$ with $1\le i_1<i_2<i_3\le \frac m2$ such that
\begin{multline} \label{eq:E8''''Civ}
\ell_1:=\ell_T(w_{i_1})=\ell_T(w_{i_1+\frac m2})\ge1,\quad
\ell_2:=\ell_T(w_{i_2})=\ell_T(w_{i_2+\frac m2})\ge1,\\
\ell_3:=\ell_T(w_{i_3})=\ell_T(w_{i_3+\frac m2})\ge1,
\quad\text{and}\quad
\ell_1+\ell_2+\ell_3\le4,
\end{multline}
and the other $w_j$'s, $1\le j\le m$, are equal to $\ep$,
\item[(v)]
there are $i_1,i_2,i_3,i_4$ with $1\le i_1<i_2<i_3<i_4\le \frac m2$ such that
\begin{multline} \label{eq:E8''''Cv}
\ell_T(w_{i_1})=\ell_T(w_{i_2})=\ell_T(w_{i_3})=\ell_T(w_{i_4})\\=
\ell_T(w_{i_1+\frac m2})=\ell_T(w_{i_2+\frac m2})=
\ell_T(w_{i_3+\frac m2})=\ell_T(w_{i_4+\frac m2})=1,
\end{multline}
and all other $w_j$'s are equal to $\ep$.
\end{enumerate}
\reseteqn

Moreover, since $(w_0;w_1,\dots,w_m)\in NC^m(E_8)$, we must have 
$w_iw_{i+\frac {m}2}\le_T c$,
respectively 
$w_{i_1}w_{i_2}w_{i_1+\frac {m}2}w_{i_2+\frac {m}2}\le_T c$,
respectively 
$$w_{i_1}w_{i_2}w_{i_3}
w_{i_1+\frac {m}2}w_{i_2+\frac {m}2}w_{i_3+\frac {m}2}\le_T c,$$
respectively 
$$w_{i_1}w_{i_2}w_{i_3}w_{i_4}
w_{i_1+\frac {m}2}w_{i_2+\frac {m}2}w_{i_3+\frac {m}2}
w_{i_4+\frac {m}2}=c.$$
Together with equations~\eqref{eq:E8''''A}--\eqref{eq:E8''''C}, 
this implies that
\begin{equation} \label{eq:E8''''D}
w_i=c^{15}w_ic^{-15}\quad\text{and}\quad 
w_i(c^7w_ic^{-7})\le_T c,
\end{equation}
respectively that
\begin{equation} \label{eq:E8''''E}
w_{i_1}=c^{15}w_{i_1}c^{-15},\quad 
w_{i_2}=c^{15}w_{i_2}c^{-15},\quad\text{and}\quad  w_{i_1}w_{i_2}(c^7w_{i_1}c^{-7})(c^7w_{i_2}c^{-7})\le_T c,
\end{equation}
respectively that
\begin{multline} \label{eq:E8''''F}
w_{i_1}=c^{15}w_{i_1}c^{-15},\quad 
w_{i_2}=c^{15}w_{i_2}c^{-15},\quad 
w_{i_3}=c^{15}w_{i_3}c^{-15},\\
\quad\text{and}\quad w_{i_1}w_{i_2}w_{i_3}
(c^7w_{i_1}c^{-7})(c^7w_{i_2}c^{-7})(c^7w_{i_3}c^{-7})\le_T c,
\end{multline}
respectively that
\begin{multline} \label{eq:E8''''G}
w_{i_1}=c^{15}w_{i_1}c^{-15},\quad 
w_{i_2}=c^{15}w_{i_2}c^{-15},\quad 
w_{i_3}=c^{15}w_{i_3}c^{-15},\quad 
w_{i_4}=c^{15}w_{i_4}c^{-15},\\
\quad\text{and}\quad w_{i_1}w_{i_2}w_{i_3}w_{i_4}
(c^7w_{i_1}c^{-7})(c^7w_{i_2}c^{-7})(c^7w_{i_3}c^{-7})
(c^7w_{i_4}c^{-7})=c.
\end{multline}
Here, the first equation in \eqref{eq:E8''''D},
the first two in \eqref{eq:E8''''E},
the first three in \eqref{eq:E8''''F},
and the first four in \eqref{eq:E8''''G}, are all automatically 
satisfied due to Lemma~\ref{lem:4} with $d=2$.

With the help of Stembridge's {\sl Maple} package {\tt coxeter}
\cite{StemAZ}, one obtains

\begin{enumerate}
\item[---] 45 solutions for $w_i$ in 
\eqref{eq:E8''''D} with $\ell_T(w_i)=1$,
\item[---]  150 solutions for $w_i$ in 
\eqref{eq:E8''''D} with $\ell_T(w_i)=2$ and $w_i$ of type $A_1^2$,
\item[---]  100 solutions for $w_i$ in 
\eqref{eq:E8''''D} with $\ell_T(w_i)=2$ and $w_i$ of type $A_2$,
\item[---]  75 solutions for $w_i$ in 
\eqref{eq:E8''''D} with $\ell_T(w_i)=3$ and $w_i$ of type $A_1^3$,
\item[---]  165 solutions for $w_i$ in 
\eqref{eq:E8''''D} with $\ell_T(w_i)=3$ and $w_i$ of type $A_1*A_2$,
\item[---]  90 solutions for $w_i$ in 
\eqref{eq:E8''''D} with $\ell_T(w_i)=3$ and $w_i$ of type $A_3$,
\item[---]  15 solutions for $w_i$ in 
\eqref{eq:E8''''D} with $\ell_T(w_i)=4$ and $w_i$ of type 
$A_1^2*A_2$,
\item[---]  45 solutions for $w_i$ in 
\eqref{eq:E8''''D} with $\ell_T(w_i)=4$ and $w_i$ of type $A_1*A_3$;
\item[---]  5 solutions for $w_i$ in 
\eqref{eq:E8''''D} with $\ell_T(w_i)=4$ and $w_i$ of type $A_2^2$,
\item[---]  18 solutions for $w_i$ in 
\eqref{eq:E8''''D} with $\ell_T(w_i)=4$ and $w_i$ of type $A_4$,
\item[---]  5 solutions for $w_i$ in 
\eqref{eq:E8''''D} with $\ell_T(w_i)=4$ and $w_i$ of type $D_4$.
\end{enumerate}

\noindent
Each of them gives rise to $m/2$ elements of
$\Fix_{NC^m(E_8)}(\phi^{p})$ since $i$ ranges from $1$ to $m/2$.
There are no solutions  for $w_i$ in \eqref{eq:E8''''D} with 
$w_i$ of type $A_1^4$.

Letting the computer find all solutions in cases (iii)--(v) would
take years. However, the number of these solutions can be computed
from our knowledge of the solutions in Case~(ii) according to type,
if this information is combined with the decomposition numbers in
the sense of \cite{KratCB,KratCF,KrMuAB} (see the end of
Section~\ref{sec:prel}) and some elementary 
(multiset) permutation counting. The decomposition numbers for
$A_2$, $A_3$, $A_4$, and $D_4$ of which we make use can be found in
the appendix of \cite{KratCF}.

To begin with, the number of
solutions of \eqref{eq:E8''''E} with $\ell_1=\ell_2=1$ is equal to
$$
n_{1,1}:=150\cdot 2+100\cdot N_{A_2}(A_1,A_1)=600,
$$
since an element of type $A_1^2$ can be decomposed in two ways 
into a product of two elements of absolute length $1$, while for
an element of type $A_2$ this can be done in $ N_{A_2}(A_1,A_1)=3$
ways. Similarly, the number of
solutions of \eqref{eq:E8''''E} with $\ell_1=2$ and $\ell_2=1$ is equal to
$$
n_{2,1}:=75\cdot 3+165\cdot(1+N_{A_2}(A_1,A_1))
+90\cdot N_{A_3}(A_2,A_1)=1425,
$$
the number of
solutions of \eqref{eq:E8''''E} with $\ell_1=3$ and $\ell_2=1$ is equal to
\begin{multline*}
n_{3,1}:=15\cdot (2+N_{A_2}(A_1,A_1))
+45\cdot(1+N_{A_3}(A_2,A_1))
+5\cdot(2N_{A_2}(A_1,A_1))\\
+18\cdot(N_{A_4}(A_3,A_1)+N_{A_4}(A_1*A_2,A_1))
+5\cdot(N_{D_4}(A_3,A_1)+N_{D_4}(A_1^3,A_1))=660,
\end{multline*}
the number of
solutions of \eqref{eq:E8''''E} with $\ell_1=\ell_2=2$ is equal to
\begin{multline*}
n_{2,2}:=15\cdot (2+2N_{A_2}(A_1,A_1))
+45\cdot(2N_{A_3}(A_2,A_1))
+5\cdot(2+N_{A_2}(A_1,A_1)^2)\\
+18\cdot(N_{A_4}(A_2,A_2)+N_{A_4}(A_1^2,A_1^2)+2N_{A_4}(A_2,A_1^2))\\
+5\cdot(N_{D_4}(A_2,A_2)+2N_{D_4}(A_2,A_1^2))=1195,
\end{multline*}
the number of
solutions of \eqref{eq:E8''''F} with $\ell_1=\ell_2=\ell_3=1$ is 
equal to
$$
n_{1,1,1}:=75\cdot 3!+165\cdot(3N_{A_2}(A_1,A_1))
+90N_{A_3}(A_1,A_1,A_1)=3375,
$$
the number of
solutions of \eqref{eq:E8''''F} with $\ell_1=2$ and 
$\ell_2=\ell_3=1$ is equal to
\begin{multline*}
n_{2,1,1}:=15\cdot (2+N_{A_2}(A_1,A_1)+2\cdot2\cdot N_{A_2}(A_1,A_1))
+45\cdot(2N_{A_3}(A_2,A_1)+N_{A_3}(A_1,A_1,A_1))\\
+5\cdot(2N_{A_2}(A_1,A_1)+2N_{A_2}(A_1,A_1)^2)
+18\cdot(N_{A_4}(A_2,A_1,A_1)+N_{A_4}(A_1^2,A_1,A_1))\\
+5\cdot(N_{D_4}(A_2,A_1,A_1)+N_{D_4}(A_1^2,A_1,A_1))=2850,
\end{multline*}
and the number of
solutions of \eqref{eq:E8''''G} is equal to
\begin{multline*}
n_{1,1,1,1}:=15\cdot (12N_{A_2}(A_1,A_1))
+45\cdot(4N_{A_3}(A_1,A_1,A_1))
+5\cdot(6N_{A_2}(A_1,A_1)^2)\\
+18\cdot N_{A_4}(A_1,A_1,A_1,A_1)
+5\cdot N_{D_4}(A_1,A_1,A_1,A_1)=6750.
\end{multline*}

In total, we obtain 
\begin{multline*}
1+(45+150+100+75+165+90+15+45+5+18+5)\frac m2
+(n_{1,1}+2n_{2,1}+2n_{3,1}+n_{2,2})\binom {m/2}2\\
+(n_{1,1,1}+3n_{2,1,1})\binom {m/2}3
+n_{1,1,1,1}\binom {m/2}4
=\frac {(5m+4)(3m+2)(5m+2)(15m+4)}{64}
\end{multline*}
elements in
$\Fix_{NC^m(E_8)}(\phi^p)$, which agrees with the limit in
\eqref{eq:E8.11}.


\section{Cyclic sieving II}
\label{sec:siev2}

In this section we present the second cyclic sieving conjecture due to
Bessis and Reiner \cite[Conj.~6.5]{BeReAA}.

Let $\psi:NC^m(W)\to NC^m(W)$ be the map defined by
\begin{equation} \label{eq:psi}
(w_0;w_1,\dots,w_m)\mapsto
\big(cw_mc^{-1};w_0,w_1,\dots,w_{m-1}\big).
\end{equation}
For $m=1$, we have $w_0=cw_1^{-1}$, so that
this action reduces to the inverse of the {\it Kreweras complement\/} 
$K_{\text {id}}^{c}$ as defined by Armstrong
\cite[Def.~2.5.3]{ArmDAA}.

It is easy to see that $\psi^{(m+1)h}$ acts as the identity,
where $h$ is the Coxeter number of $W$ (see \eqref{eq2:Aktion} below).
By slight abuse of notation as before, let $C_2$ be the cyclic group of order
$(m+1)h$ generated by $\psi$. 

Given these definitions, we are now in the position to state
the second cyclic sieving conjecture of Bessis and Reiner.
By the results of \cite{KratCG} and of this paper, it becomes the
following theorem.

\begin{theorem} \label{thm2:1}
For an irreducible well-generated complex reflection group $W$ and any $m\ge1$, 
the triple $(NC^m(W),\Cat^m(W;q),C_2)$, where $\Cat^m(W;q)$ is
the $q$-analogue of the Fu\ss--Catalan number defined in
\eqref{eq:FCZahl}, exhibits the cyclic sieving phenomenon.
\end{theorem}

By definition of the cyclic sieving phenomenon, we have to
prove that 
\begin{equation} \label{eq2:1}
\vert\Fix_{NC^m(W)}(\psi^{p})\vert = 
\Cat^m(W;q)\big\vert_{q=e^{2\pi i p/(m+1)h}}, 
\end{equation}
for all $p$ in the range $0\le p<(m+1)h$.

\section{Auxiliary results II}
\label{sec:aux2}

This section collects several auxiliary results which allow us to
reduce the problem of proving Theorem~\ref{thm2:1}, respectively the
equivalent statement
\eqref{eq2:1}, for the 26 exceptional groups listed in
Section~\ref{sec:prel} to a finite problem. 
The corresponding lemmas, Lemmas~\ref{lem2:1}--\ref{lem2:8},
are analogues of Lemmas~\ref{lem:1}--\ref{lem:3} and
\ref{lem:6}--\ref{lem:8} in Section~\ref{sec:aux1}.

\medskip
Let $p=a(m+1)+b$, $0\le b<m+1$. We have
\begin{multline}
\psi^p\big((w_0;w_1,\dots,w_m)\big)\\
=(c^{a+1}w_{m-b+1}c^{-a-1};c^{a+1}w_{m-b+2}c^{-a-1},
\dots,c^{a+1}w_{m}c^{-a-1},\\
c^{a}w_{0}c^{-a},\dots,
c^{a}w_{m-b}c^{-a}\big).
\label{eq2:Aktion}
\end{multline}

\begin{lemma} \label{lem2:1}
It suffices to check \eqref{eq2:1} for $p$ a divisor of $(m+1)h$.
More precisely, let $p$ be a divisor of $(m+1)h$, and let $k$ be another positive integer with 
$\gcd(k,(m+1)h/p)=1$, then we have
\begin{equation} \label{eq2:2}
\Cat^m(W;q)\big\vert_{q=e^{2\pi i p/(m+1)h}}
= 
\Cat^m(W;q)\big\vert_{q=e^{2\pi i kp/(m+1)h}}
\end{equation}
and
\begin{equation} \label{eq2:3}
\vert\Fix_{NC^m(W)}(\psi^{p})\vert =
\vert\Fix_{NC^m(W)}(\psi^{kp})\vert 	.
\end{equation}
\end{lemma}

\begin{proof}
For \eqref{eq2:3}, this follows in the same way as \eqref{eq:3}
in Lemma~\ref{lem:1}.

For \eqref{eq2:2}, we must argue differently than in 
Lemma~\ref{lem:1}. Let us write $\zeta=e^{2\pi i p/(m+1)h}$.
For a given group $W$, we write $S_1(W)$ for the set of all
indices $i$ such that $\zeta^{d_i-h}=1$, and we write
$S_2(W)$ for the set of all
indices $i$ such that $\zeta^{d_i}=1$. 
By the rule of de l'Hospital, we have
\begin{equation} \label{eq:S1S2}
\Cat^m(W;q)\big\vert_{q=e^{2\pi i p/(m+1)h}}=
\begin{cases}
0&\text{if }\vert S_1(W)\vert>\vert S_2(W)\vert,\\
\frac {\prod_{i\in S_1(W)}(mh+d_i)} 
{\prod_{i\in S_2(W)} d_i}
\frac {\prod_{i\notin S_1(W)}(1-\zeta^{d_i-h})} 
{\prod_{i\notin S_2(W)}(1-\zeta^{d_i})},
&\text{if }\vert S_1(W)\vert=\vert S_2(W)\vert.
\end{cases}
\end{equation}
Since, by Theorem~\ref{thm:0}, $\Cat^m(W;q)$ is a polynomial in $q$,
the case 
$\vert S_1(W)\vert<\vert S_2(W)\vert$ cannot occur.

We claim 
that, for the case where $\vert S_1(W)\vert=\vert S_2(W)\vert$, 
the factors in the quotient of products
$$
\frac {\prod_{i\notin S_1(W)}(1-\zeta^{d_i-h})} 
{\prod_{i\notin S_2(W)}(1-\zeta^{d_i})}
$$
cancel pairwise. If we assume the correctness of the claim,
it is obvious that we get the same result
if we replace $\zeta$ by $\zeta^k$, where $\gcd(k,(m+1)h/p)=1$,
hence establishing \eqref{eq2:2}.

In order to see that our claim is indeed valid,
we proceed in a case-by-case fashion,
making appeal to the classification of irreducible well-generated complex 
reflection groups, which we recalled in Section~\ref{sec:prel}.
First of all, since $d_n=h$, the set $S_1(W)$ is always non-empty
as it contains the element $n$. Hence, if we want to have
$\vert S_1(W)\vert=\vert S_2(W)\vert$, the set $S_2(W)$ must be
non-empty as well. In other words, the integer $(m+1)h/p$ must divide at least
one of the degrees $d_1,d_2,\dots,d_n$. In particular, this implies that,
for each fixed reflection group $W$ of exceptional type, only a finite number
of values of $(m+1)h/p$ has to be checked. Writing $M$ for $(m+1)h/p$,
what needs to be checked is whether the 
{\it multi}sets (that is, multiplicities of elements must be taken
into account)
$$
\{(d_i-h)\text{ mod }M:i\notin S_1(W)\}
\quad 
\text{and}
\quad 
\{d_i\text{ mod }M:i\notin S_2(W)\}
$$
are the same. Since, for a fixed irreducible well-generated complex
reflection group, there is only a finite number of possibilities for
$M$, this amounts to a routine verification.
\end{proof}

\begin{lemma} \label{lem2:2}
Let $p$ be a divisor of $(m+1)h$.
If $p$ is divisible by $m+1$, then \eqref{eq2:1} is true.
\end{lemma}

We leave the proof to the reader as it is completely analogous 
to the proof of Lemma~\ref{lem:2}.
%

\begin{lemma} \label{lem2:3}
Equation \eqref{eq2:1} holds for all divisors $p$ of $m+1$.
\end{lemma}

\begin{proof}
We have
$$
\Cat^m(W;q)\big\vert_{q=e^{2\pi i p/(m+1)h}}=
\begin{cases}
0&\text{if }p<m+1,\\
m+1&\text{if }p=m+1.
\end{cases}
$$
Here, the first case follows from \eqref{eq:S1S2} and the fact
that we have $S_1(W)\supseteq\{n\}$ 
and $S_2(W)=\emptyset$ if $p\mid (m+1)$ and
$p<m+1$.

On the other hand, if $(w_0;w_1,\dots,w_m)$ is fixed by $\psi^p$,
then one can apply an argument similar to that in Lemma~\ref{lem:3}
with any $w_i$ taking the role of $w_1$, $0\le i\le m$. It follows 
that if
$p=m+1$, the set $\Fix_{NC^m(W)}(\psi^p)$ consists of the 
$m+1$ elements $(w_0;w_1,\dots,w_m)$ obtained by choosing 
$w_i=c$ for a particular $i$ between $0$ and $m$, all other $w_j$'s 
being equal to $\ep$. If 
$p<m+1$, then
there is no element in $\Fix_{NC^m(W)}(\psi^p)$.
\end{proof}

\begin{lemma} \label{lem2:6}
Let $W$ be an irreducible 
well-generated complex reflection group of rank $n$, 
and let $p=m_1h_1$ be a divisor of $(m+1)h$, where $m+1=m_1m_2$ and
$h=h_1h_2$. We assume that 
$\gcd(h_1,m_2)=1$. Suppose that
Theorem~{\em\ref{thm2:1}} has already been verified for all 
irreducible well-generated
complex reflection groups with rank $<n$. 
If $h_2$ does not divide all degrees $d_i$,
then equation~\eqref{eq2:1} is satisfied.
\end{lemma}

We leave the proof to the reader as it is completely analogous 
to the proof of Lemma~\ref{lem:6}.

\begin{lemma} \label{lem2:7}
Let $W$ be an irreducible 
well-generated complex reflection group of rank $n$, 
and let $p=m_1h_1$ be a divisor of $(m+1)h$, where $m+1=m_1m_2$ and
$h=h_1h_2$. We assume that 
$\gcd(h_1,m_2)=1$. If $m_2>n$ then
$$
\Fix_{NC^m(W)}(\psi^{p})=\emptyset.
$$
\end{lemma}

We leave the proof to the reader as it is analogous 
to the proof of Lemma~\ref{lem:7}.


\begin{remark} \label{rem2:1}
By applying the same reasoning as in Remark~\ref{rem:1} with Lemmas~\ref{lem:6}
and \ref{lem:7} replaced by Lemmas~\ref{lem2:6} and \ref{lem2:7}, respectively,
it follows that we only need to check \eqref{eq2:1} for pairs $(m_2,h_2)$
satisfying \eqref{eq:restr} and $m_2\le n$. This reduces the problem
to a finite number of choices.
%
%
%
\end{remark}

\begin{lemma} \label{lem2:8}
Let $W$ be an irreducible 
well-generated complex reflection group of rank $n$ 
with the property that $d_i\mid h$ for $i=1,2,\dots,n$.
Then Theorem~{\em\ref{thm2:1}} is true for this group $W$.
\end{lemma}

\begin{proof}

%
%
Proceeding in a fashion analogous to the beginning of the proof of
Lemma~\ref{lem:8}, we may restrict to the case where $p\mid (m+1)h$ and
$(m+1)h/p$ does not divide any of the $d_i$'s.
In this case, it follows from \eqref{eq:S1S2} and the fact that
we have  $S_1(W)\supseteq\{n\}$ 
and $S_2(W)=\emptyset$ that
the right-hand side of \eqref{eq2:1} equals $0$. 
Inspection of the classification of all irreducible 
well-generated complex reflection groups, which we recalled in
Section~\ref{sec:prel}, reveals that all groups
satisfying the hypotheses of the lemma have rank $n\le2$.
Except for the groups contained in the infinite series $G(d,1,n)$ 
and $G(e,e,n)$ for which Theorem~\ref{thm:1} has been established in
\cite{KratCG}, these are the groups $G_5,G_6,G_9,G_{10},
G_{14},G_{17},G_{18},G_{21}$. 
The verification of \eqref{eq2:1} can be done in a similar fashion as
in the proof of Lemma~\ref{lem:8}. We illustrate this by going through 
the case of the group $G_6$. In analogy with the earlier situation, 
we note that Lemma~\ref{lem2:7}
implies that equation~\eqref{eq2:1} holds if $m_2>2$, so that
in the following arguments we may assume that $m_2=2$.

\smallskip
{\sc Case $G_6$}. The degrees are $4,12$, and
therefore, according to Remark~\ref{rem2:1}, we need only consider
the case where $h_2=4$ and $m_2=2$, that is, $p=3(m+1)/2$. 
Then the action of $\psi^p$ is given by
\begin{equation} \label{eq2:3m2Aktion}
\psi^p\big((w_0;w_1,\dots,w_m)\big)
=(c^{2}w_{\frac {m+1}2}c^{-2};
c^{2}w_{\frac {m+3}2}c^{-2},
\dots,c^{2}w_{m}c^{-2},
cw_{0}c^{-1},\dots,
cw_{\frac {m-1}2}c^{-1}\big).
\end{equation}
If $(w_0;w_1,\dots,w_m)$ is fixed by $\psi^p$, 
there must exist an $i$
with $0\le i\le \frac {m-1}2$ such that $\ell_T(w_i)=1$,
$w_icw_{i}c^{-1}=c$, and all $w_j$, $j\ne i,\frac {m+1}2+i$,
equal $\ep$. However, 
with the help of {\tt CHEVIE}, one verifies that there is no such
solution to this equation. Hence,
the left-hand side of \eqref{eq2:1} is equal to $0$, as required.

\smallskip
This completes the proof of the lemma.
\end{proof}

\section{Exemplification of case-by-case verification of Theorem~\ref{thm2:1}}
\label{sec:Beweis2}

It remains to verify Theorem~\ref{thm2:1} for the groups
$G_4,G_8,G_{16},G_{20},
G_{23}=H_3,G_{24},G_{25},\break G_{26},G_{27},
G_{28}=F_4,G_{29},G_{30}=H_4,G_{32},
G_{33},G_{34},G_{35}=E_6,G_{36}=E_7,G_{37}=E_8$.
All details can be found in \cite[Sec.~9]{KrMuAE}.
We content ourselves with discussing the case of the group $G_{24}$,
as this suffices to convey the flavour of
the necessary computations.

In order to simplify our considerations,
it should be observed that the action of $\psi$ (given in 
\eqref{eq:psi}) is exactly the same as the action of $\phi$
(given in \eqref{eq:phi}) with $m$ replaced by $m+1$ 
{\it on the components $w_1,w_2,\dots,w_{m+1}$}, that is,
if we disregard the $0$-th component of the elements of the
generalised non-crossing partitions involved. The only difference
which arises is that, while the $(m+1)$-tuples 
$(w_0;w_1,\dots,w_m)$ in \eqref{eq:psi} must satisfy 
$w_0w_1\cdots w_m=c$, for  $w_1,w_2,\dots,w_{m+1}$ in \eqref{eq:phi}
we only must have $w_1w_2\cdots w_{m+1}\le_T c$. 
Consequently, we may use the counting results from 
Section~\ref{sec:Beweis1}, except that we have to restrict our
attention to those elements $(w_0;w_1,\dots,w_m,w_{m+1})\in 
NC^{m+1}(W)$ for which $w_1w_2\cdots w_{m+1}=c$, or, equivalently,
$w_0=\ep$.


\subsection*{\sc Case $G_{24}$}
The degrees are $4,6,14$, and hence we have
$$
\Cat^m(G_{24};q)=\frac 
{[14m+14]_q\, [14m+6]_q\, [14m+4]_q} 
{[14]_q\, [6]_q\, [4]_q} .
$$
Let $\zeta$ be a $14(m+1)$-th root of unity. 
The following cases on the right-hand side of \eqref{eq2:1}
occur:
{\refstepcounter{equation}\label{eq2:G24}}
\alphaeqn
\begin{align} 
\label{eq2:G24.2}
\lim_{q\to\zeta}\Cat^m(G_{24};q)&=m+1,
\quad\text{if }\zeta=\zeta_{14},\zeta_7,\\
\label{eq2:G24.3}
\lim_{q\to\zeta}\Cat^m(G_{24};q)&=\tfrac {7m+7}3,
\quad\text{if }\zeta=\zeta_{6},\zeta_3,\ 3\mid (m+1),\\
\label{eq2:G24.5}
\lim_{q\to\zeta}\Cat^m(G_{24};q)&=\Cat^m(G_{24}),
\quad\text{if }\zeta=-1\text{ or }\zeta=1,\\
\label{eq2:G24.1}
\lim_{q\to\zeta}\Cat^m(G_{24};q)&=0,
\quad\text{otherwise.}
\end{align}
\reseteqn

We must now prove that the left-hand side of \eqref{eq2:1} in
each case agrees with the values exhibited in 
\eqref{eq2:G24}. The only cases not covered by
Lemma~\ref{lem2:2} are the ones in \eqref{eq2:G24.3}
and \eqref{eq2:G24.1}. 
On the other hand, the only cases left
to consider according to Remark~\ref{rem2:1} are 
the cases where $h_2=1$ and $m_2=3$, $h_2=2$ and $m_2=3$, 
and $h_2=m_2=2$. 
These correspond to the choices $p=14(m+1)/3$,
$p=7(m+1)/3$, respectively $p=7(m+1)/2$.
The first two cases belong to \eqref{eq2:G24.3}, while
$p=7(m+1)/2$ belongs to \eqref{eq2:G24.1}.

In the case that $p=7(m+1)/3$, 
the action of $\psi^p$ is given by
\begin{multline*}
\psi^p\big((w_0;w_1,\dots,w_m)\big)\\
=(c^{3}w_{\frac {2m+2}3}c^{-3};c^{3}w_{\frac {2m+5}3}c^{-3},
\dots,c^{3}w_{m}c^{-3},
c^2w_{0}c^{-2},\dots,
c^2w_{\frac {2m-1}3}c^{-2}\big).
\end{multline*}
Hence, for an $i$ with $0\le i\le \frac {m-2} {3}$,
we must find an element $w_i=t_1$, where $t_1$ satisfies
\eqref{eq:G24D}, so that we can set $w_{i+\frac {m+1} {3}}
=c^2t_1c^{-2}$, $w_{i+\frac {2m+2} {3}}
=c^4t_1c^{-4}$, and all other $w_j$'s equal to $\ep$.
We have found
seven solutions to the counting problem \eqref{eq:G24D},
and each of them gives rise to $(m+1)/3$ elements in
$\Fix_{NC^m(G_{24})}(\psi^p)$ 
since the index $i$ ranges from $0$ to $(m-2)/3$. 

On the other hand, if $p=14(m+1)/3$, then the action of $\psi^p$ is given by
\begin{multline*}
\psi^p\big((w_0;w_1,\dots,w_m)\big)\\
=(c^{5}w_{\frac {m+1}3}c^{-5};c^{5}w_{\frac {m+4}3}c^{-5},
\dots,c^{5}w_{m}c^{-5},
c^4w_{0}c^{-4},\dots,
c^4w_{\frac {m-2}3}c^{-4}\big).
\end{multline*}
By Lemma~\ref{lem:4}, every element of $NC(W)$ is fixed under 
conjugation by $c^7$, and, thus, 
the equations for $t_1$ in this case are the same as in
the previous one where $p=7(m+1)/3$.

Hence, in either case, we obtain 
$7\frac {m+1}3=\frac {7m+7}3$ elements in
$\Fix_{NC^m(G_{24})}(\psi^p)$, which agrees with the limit in
\eqref{eq2:G24.3}.

If $p=7(m+1)/2$, the relevant counting problem is \eqref{eq:G24''D}.
However, no element
$(w_0;w_1,\dots,w_m)\in \Fix_{NC^m(G_{24})}(\psi^p)$ can be
produced in this way since the counting problem imposes the
restriction that $\ell_T(w_0)+\ell_T(w_1)+\dots+\ell_T(w_m)$
be even, which contradicts the fact that $\ell_T(c)=n=3$. 
This is in agreement with the limit in
\eqref{eq2:G24.1}.

\section*{Acknowledgements}
The authors thank an anonymous referee for a very careful 
reading of the original manuscript, and for numerous
pertinent suggestions which have
helped to considerably improve the original manuscript.

\end{document}
